\documentclass[camera,extra]{gji}
\usepackage[dvips]{epsfig,color,graphicx} 
\usepackage{times} 
\usepackage{makeidx} 
\usepackage{url,phonetic} 
\usepackage{fancyhdr} 
\usepackage{letterspace}
\usepackage{eurosym} 
\usepackage{subfigure}
\usepackage{mathrsfs} 
\usepackage[fleqn]{amsmath} 
\title[Slepian functions and the polar gap]{Spherical Slepian functions and the polar gap in geodesy}
\author[F.~J. Simons and F.~A.~Dahlen]{Frederik J.~Simons$^1$ and F.~A.~Dahlen$^2$\\$^1$ Department of Earth Sciences, University College London, Gower Street, London WC1E 6BT, U.~K.\\$^2$ Department of Geosciences, Princeton University, Guyot Hall, Princeton, NJ 08544, U.~S.~A.} 
\date{Accepted \today. Received \today; in original form \today}
\pagerange{\pageref{firstpage}--\pageref{lastpage}}
\volume{XXX}
\pubyear{XXXX}
\newcommand{\domg}{\,d\Omega} 
\newcommand{\fracd}[2]{\frac{\displaystyle{#1}}{\displaystyle{#2}}} 
\newcommand{\intr}{\int_R} 
\newcommand{\intbr}{\int_{\bar{R}}} 
\newcommand{\into}{\int_\Omega} 
\newcommand{\covsl}{\langle s_{lm}s_{l'm'}\rangle} %
\newcommand{\covsr}{\langle s(\rhat)s(\rhat')\rangle} %
\newcommand{\covsa}{\langle s_\alpha s_\beta\rangle} %
\newcommand{\sfgup}{\sfg_{\uparrow}} %
\newcommand{\sfgdn}{\sfg_{\downarrow}} %
\newcommand{\gup}{g_{\uparrow}} %
\newcommand{\gdn}{g_{\downarrow}} %
\newcommand{\gdnlm}{g_{\downarrow\hspace*{0.1em}lm}} %
\newcommand{\gdna}{g_{\downarrow\hspace*{0.1em}\alpha}} %
\newcommand{\gupa}{g_{\uparrow\hspace*{0.1em}\alpha}} %
\newcommand{\siup}{s_{\uparrow}}%
\newcommand{\siuph}{\shat_{\uparrow}}%
\newcommand{\siupa}{s_{\uparrow\hspace*{0.05em}\alpha}}%
\newcommand{\hsiupa}{\hat{s}_{\uparrow\hspace*{0.05em}\alpha}}%
\newcommand{\sidna}{s_{\downarrow\hspace*{0.05em}\alpha}}%
\newcommand{\sfgupa}{\sfg_{\uparrow\hspace*{0.05em}\alpha}} %
\newcommand{\sfgupb}{\sfg_{\uparrow\hspace*{0.05em}\beta}} %
\newcommand{\sfgdna}{\sfg_{\downarrow\hspace*{0.05em}\alpha}} %
\newcommand{\sfgdnb}{\sfg_{\downarrow\hspace*{0.05em}\beta}} %
\newcommand{\guplm}{g_{\uparrow\hspace*{0.1em}lm}} %
\newcommand{\suplm}{s_{\uparrow\hspace*{0.1em}lm}} %
\newcommand{\suplmh}{\hat{s}_{\uparrow\hspace*{0.1em}lm}} %
\newcommand{\suplmp}{s_{\uparrow\hspace*{0.1em}l'm'}} %
\newcommand{\suplmpp}{s_{\uparrow\hspace*{0.1em}l''m''}} %
\newcommand{\lstart}{[\lambda_\alpha+\eta(1-\lambda_\alpha)]^{-1}} %
\newcommand{\lstar}{\lambda_\alpha^*(\eta)} %
\newcommand{\lstarR}{\lambda_\alpha^*(\eta_R)} %
\newcommand{\lstarO}{\lambda_\alpha^*(\eta_\Omega)} %
\newcommand{\also}{\quad\mbox{and}\quad} %
\newcommand{\dab}{\delta_{\alpha\beta}} %
\newcommand{\dlmp}{\delta_{ll'}\delta_{mm'}} %
\newcommand{\dllp}{\delta_{ll'}} %
\newcommand{\dmmp}{\delta_{mm'}} %
\newcommand{\Ylm}{Y_{lm}}  
\newcommand{\Ylmrh}{Y_{lm}(\rhat)}  %
\newcommand{\Ylmrhp}{Y_{lm}(\rhat')}  %
\newcommand{\Ylmprh}{Y_{l'm'}(\rhat)}  %
\newcommand{\Ylmp}{Y_{l'm'}} 
\newcommand{\Xlm}{X_{lm}}  %
\newcommand{\Xlmth}{X_{lm}(\theta)}  %
\newcommand{\Xlpm}{X_{l'm}} %
\newcommand{\Dlmlmp}{D_{lm,l'm'}} %
\newcommand{\bDlmlmp}{\bar{D}_{lm,l'm'}} %
\newcommand{\Dlmpplmp}{{D}_{l''m'',l'm'}} %
\newcommand{\Dlmplmpp}{{D}_{l'm',l''m''}} %
\newcommand{\Dlmlmpp}{{D}_{lm,l''m''}} %
\newcommand{\Dllpm}{D^m_{ll'}} %
\newcommand{\bDllp}{\bar{D}_{ll'}} %
\newcommand{\Dllppm}{D^m_{ll''}} %
\newcommand{\Dllpp}{D_{ll''}} %
\newcommand{\Dlpplpm}{D^m_{l''l'}} %
\newcommand{\Dlpplp}{D_{l''l'}} %
\newcommand{\Almlmp}{A_{lm,l'm'}} 
\newcommand{\Drhrhp}{D(\rhat,\rhat')} 
\newcommand{\grhp}{g(\rhat')} 
\newcommand{\garh}{g_\alpha(\rhat)} 
\newcommand{\galm}{g_{\alpha,lm}} 
\newcommand{\galmp}{g_{\alpha,l'm'}} 
\newcommand{\halm}{h_{\alpha,lm}} 
\newcommand{\glma}{g_{\alpha,lm}} 
\newcommand{\glmb}{g_{\beta,lm}} 
\newcommand{\Lpot}{(L+1)^2} 
\newcommand{\sumapot}{\sum_{\alpha=1}^{\Lpot}} 
\newcommand{\sumbpot}{\sum_{\beta=1}^{\Lpot}} 
\newcommand{\sumak}{\sum_{\alpha=1}^{J}} 
\newcommand{\sumakR}{\sum_{\alpha>J}^{\Lpot}} 
\newcommand{\barr}{\bar{R}} 
\newcommand{\suml}{\sum\limits}  
\newcommand{\sumelpp}{\suml_{l''=0}^{L}} 
\newcommand{\sumsh}{\suml_{l=0}^{\infty}\suml_{m=-l}^{l}} 
\newcommand{\sumshort}{\suml_{lm}^{\infty}} %
\newcommand{\sumshortp}{\suml_{l'm'}^{\infty}} %
\newcommand{\sumshortpp}{\suml_{l''m''}^{\infty}} %
\newcommand{\sumshL}{\suml_{l=0}^{L}\suml_{m=-l}^{l}} 
\newcommand{\sumshortL}{\suml_{lm}^{L}} 
\newcommand{\sumshortLp}{\suml_{l'm'}^{L}} 
\newcommand{\sumshortLpp}{\suml_{l''m''}^{L}} 
\newcommand{\sumshR}{\suml_{l>L}^{\infty}\suml_{m=-l}^{l}} 
\newcommand{\sumshortR}{\suml_{lm>L}^{\infty}} 
\newcommand{\sumshortRp}{\suml_{l'm'>L}^{\infty}} 
\newcommand{\tlofp}{\left(\frac{2l+1}{4\pi}\right)} 
\newcommand{\Plm}{P_{lm}}  
\newcommand{\sfD}{{\textsf D}} %
\newcommand{\sfT}{{\textsf T}} %
\newcommand{\sfg}{{\textsf{\small{g}}}} %
\newcommand{\T}{^{\sf{\sst{T}}}} %
\newcommand{\sfA}{{\textsf A}} %
\newcommand{\sfI}{{\textsf I}} %
\newcommand{\sfbD}{{\bar{\textsf D}}} %
\newcommand{\wigner}[6] 
{\left(\barray{ccc} #1 & #2 & #3\\ #4 & #5 & #6 \earray\right)} 
\newcommand{\matd}[9] 
{\left(\barray{ccc} #1 & #2 & #3\\ #4 & #5 & #6\\ #7 & #8 & #9 \earray\right)} 
\newcommand{\detd}[9] 
{\left|\barray{ccc} #1 & #2 & #3\\ #4 & #5 & #6\\ #7 & #8 & #9 \earray\right|} 
\newcommand{\matt}[4] 
{\left(\barray{cc} #1 & #2 \\ #3 & #4 \earray\right)} 
\newcommand{\dett}[4] 
{\left|\barray{cc} #1 & #2 \\ #3 & #4 \earray\right|} 
\newcommand{\ffae}{\mathcal{R}_\alpha(\eta)}  %
\newcommand{\ffan}{\mathcal{R}_\alpha(0)}  %
\newcommand{\rhat}{\mbf{\hat{r}}} 
\newcommand{\shat}{\hat{s}} 
\newcommand{\bmp}{\begin{minipage}} 
\newcommand{\emp}{\end{minipage}} 
\newcommand{\be}{\begin{equation}} 
\newcommand{\ee}{\end{equation}} 
\newcommand{\bt}{\begin{tabular}} 
\newcommand{\et}{\end{tabular}} 
\newcommand{\btstar}{\begin{tabular*}} 
\newcommand{\etstar}{\end{tabular*}} 
\newcommand{\btx}{\begin{tabularx}} 
\newcommand{\etx}{\end{tabularx}} 
\newcommand{\bn}{\begin{enumerate}} 
\newcommand{\en}{\end{enumerate}} 
\newcommand{\bd}{\begin{description}} 
\newcommand{\ed}{\end{description}} 
\newcommand{\bdc}{\begin{document}} 
\newcommand{\edc}{\end{document}} 
\newcommand{\bi}{\begin{itemize}} 
\newcommand{\ei}{\end{itemize}} 
\newcommand{\ber}{\begin{eqnarray}} 
\newcommand{\eer}{\end{eqnarray}} 
\newcommand{\barray}{\begin{array}} 
\newcommand{\earray}{\end{array}} 
\newcommand{\mbf}{\mathbf} 
\newcommand{\mrm}{\mathrm} 
\newcommand{\mtl}{\mathit} 
\newcommand{\pl}{\partial} 
\newcommand{\nb}{\nabla} 
\newcommand{\al}{\alpha} 
\newcommand{\nnr}{\nonumber} 
\newcommand{\fns}{\footnotesize} 
\newcommand{\nns}{\normalsize} 
\newcommand{\sst}{\scriptstyle} 
\newcommand{\ssts}{\scriptscriptstyle} 
\newcommand{\rar}{\rightarrow} 
\newcommand{\sct}{\section} 
\newcommand{\ssec}{\subsection} 
\newcommand{\lbl}{\label} 
\newcommand{\Rar}{\Rightarrow} 
\begin{document} 
\newcommand{\fmat}[4] 
{\left[\barray{cc} #1 & #2\\ #3 & #4 \earray\right]} 
\hyphenation{band-limited band-lim-i-ta-tion axi-sym-met-ric
non-sto-chas-tic sto-chas-tic}
\newcommand{\hsp}{\hspace*{0.1em}}
\newcommand{\hsps}{\hspace*{0.75em}}
\newcommand{\hspm}{\hspace*{-0.1em}}
\newcommand{\sump}{\sideset{}{'}\sum}
\newcommand{\sumpp}{\sump_{l=m_p}^{L_p}}
\newcommand{\sumpn}{\sump_{n=m_p}^{l}}
\newcommand{\Alma}{C_{lm}}
\newcommand{\Almb}{C_{nm}}
\newcommand{\Alm}{(2l+1)\fracd{(l-m)!}{(l+m)!}}
\newcommand{\ssT}{\mathcal{T}}
\maketitle
\begin{summary}
The estimation of potential fields such as the gravitational or
magnetic potential at the surface of a spherical planet from noisy
observations taken at an altitude over an incomplete portion of the
globe is a classic example of an ill-posed inverse problem. Here we
show that the geodetic estimation problem has deep-seated connections
to Slepian's spatiospectral localization problem on the sphere, which
amounts to finding bandlimited spherical functions whose energy is
optimally concentrated in some closed portion of the unit sphere. This
allows us to formulate an alternative solution to the traditional
damped least-squares spherical harmonic approach in geodesy, whereby
the source field is now expanded in a truncated Slepian function basis
set. We discuss the relative performance of both methods with regard
to standard statistical measures as bias, variance and mean-square
error, and pay special attention to the algorithmic efficiency of
computing the Slepian functions on the region complementary to the
axisymmetric polar gap characteristic of satellite surveys. The ease,
speed, and accuracy of this new method makes the use of spherical
Slepian functions in earth and planetary geodesy practical.
\end{summary}
\label{firstpage}

\section{I~n~t~r~o~d~u~c~t~i~o~n}

Satellites mapping out the spatial variations of the gravitational or
magnetic fields of the Earth or other planets ideally fly on polar
orbits, uniformly covering the entire globe. Thus potential fields on
the sphere are usually expressed in spherical harmonics, basis
functions with global support. For various, especially engineering,
reasons, however, inclined orbits are favorable. These leave a ``polar
gap'': an antipodal pair of axisymmetric polar caps, typically less
than 10$^\circ$ in diameter, without any data coverage. Estimation of
spherical harmonic field coefficients from an incompletely sampled
sphere is prone to error, since the spherical harmonics are not
orthogonal over the partial domain of the cut sphere.

The historically somewhat neglected geodetic polar gap problem has
been revived by, among others, \cite{Sneeuw+97}, and recently,
\cite{Albertella+99}, who constructed a new basis of so-called Slepian 
functions \cite[after][]{Slepian83} on the sphere. These bandlimited
functions are designed to have the majority of their energy optimally
concentrated inside the latitudinal belt composed of the entire globe
minus the polar gap, i.e. the region covered by satellites. Slepian
functions are orthogonal on both the entire as well as the cut sphere,
a property that can be exploited to our advantage. Here, we study the
inverse problem of retrieving a potential field on the unit sphere
from noisy and incomplete observations made at an altitude above their
source. We derive exact expressions for the estimation error due to
the traditional method of damped least-squares spherical harmonic
analysis as well as that arising from a new approach using a truncated
set of Slepian basis functions.

We cast the geodetic estimation problem in the much wider context of
spatiospectral localization, whereby bandlimited functions are
spatially concentrated to regions of arbitrary shape on the sphere
\cite[]{Wieczorek+2005,Simons+2006a}, and derive a new
semi-analytical numerical method to calculate the spherical Slepian
functions on the latitudinal belt or its complement, the double polar
cap. Our approach requires no numerical integration, and avoids the
construction of matrices other than a tridiagonal matrix whose
elements are prescribed analytically. Finding spherical harmonic
expressions for bandlimited functions concentrated to polar caps or
latitudinal belts, as in Figure\ \ref{sddiagram}, thus becomes so
effortless as to be achievable by a handful of lines of computer code,
and the problems with numerical stability known to plague alternative
approaches
\cite[]{Albertella+99,Pail+2001} are avoided altogether. 

The key to this ``magic'' lay hidden in two little-known studies
published several decades ago: the work by \cite{Gilbert+77} on doubly
orthogonal polynomials, and that on commuting differential operators
by \cite{Grunbaum+82}. It must be remembered that one of Slepian's
main discoveries \cite[see, e.g.,][]{Slepian83} was the existence of a
second-order differential operator that commutes with the
spatiospectral localization kernel concentrating to intervals on the
real line. Cast in matrix form, finding the prolate spheroidal
functions amounts to the diagonalization of a simple tridiagonal
matrix \cite[see, e.g.][]{Percival+93}. In their study,
\cite{Gilbert+77} presented two additional commuting differential
operators, which are applicable to the concentration of Legendre
polynomials to one-and two-sided domains. \cite{Grunbaum+82} proved
that the matrix accompanying the localization to the single polar cap
is, once again, tridiagonal. Here, we show this is also the case for
the antipodal double polar cap and its complement, the latitudinal
belt. The tridiagonal matrix elements coding for the single polar cap,
and their solutions, were published by us elsewhere
\cite[]{Simons+2006a}. The expressions applicable to the double polar
cap appear here for the first time. 

The problems posed and solved in this paper are not limited to geodesy
and observations made from a satellite. In geomagnetism, our
observation level may be the Earth's surface, and the source level at
or near the core-mantle boundary. In cosmology, the unit sphere
constituting the sky is observed from the inside out, and the galactic
plane masking spacecraft measurements has the shape of a latitudinal
belt \cite[]{Tegmark96a,Hinshaw+2003}. Ground-based astronomical
measurements may be confined to a small circular patch of the sky
\cite[]{Peebles73,Tegmark95}. Finally, in planetary science, 
knowledge of the estimation statistics of properties observed over
mere portions of the planetary surface is important in the absence of
groundtruthing observations.

\section{S~t~a~t~e~m~e~n~t{\hsps}o~f{\hsps}t~h~e{\hsps}p~r~o~b~l~e~m}
\label{ps}

\begin{figure}\center
\rotatebox{0}{
\iftwocol
{\includegraphics[width=\columnwidth]{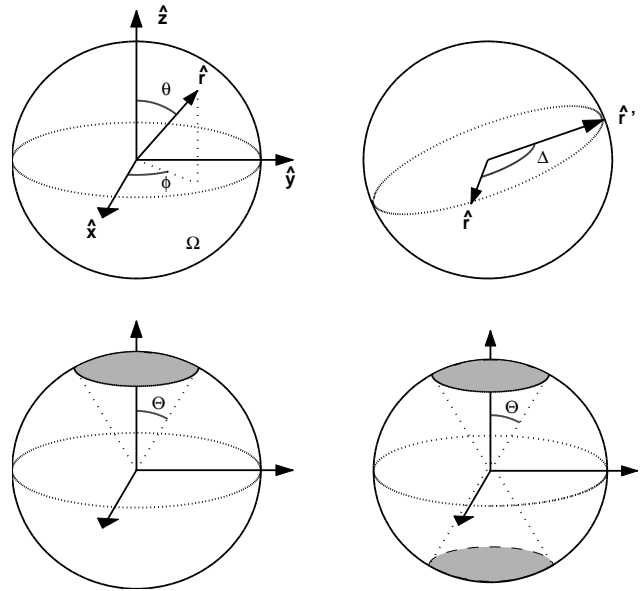}}
{\includegraphics[width=0.65\columnwidth]{sddiagram.eps}}
}
\caption{Geometry of the geodetic estimation problem, and some
symbols used in this paper. In the lower left, an axisymmetric polar
cap, shaded, of colatitudinal radius $\Theta$. In the lower right, an
antipodal pair of polar caps, shaded, representing the geodetic polar
gap.}
\label{sddiagram}  
\end{figure}

We are concerned with estimating source-level potential fields from
noise-contaminated satellite observations at an altitude over an
incomplete portion of the unit sphere. The geometry of this problem is
illustrated in Figure~\ref{sddiagram}. The unit sphere $\Omega$ on
which the unknown signal is defined is parameterized as usual in terms
of spherical coordinates, colatitude $\theta$ and longitude
$\phi$. The angular distance between two position coordinates
$\rhat=(\theta,\phi)$ and $\rhat'=(\theta',\phi')$ is denoted by
$\Delta$. In the lower right, the domain over which satellite
observations are available is left unshaded, whereas the area in which
measurements are missing is shaded grey. We denote the white region
covered by satellite tracks by $R$, and the shaded, uncovered region
by $\barr$. Although our treatment will start out quite general,
without restrictions on the shape of $R$ or $\barr$, as long as they
are complementary closed regions on the surface of the unit sphere,
the lower right panel of Figure~\ref{sddiagram} illustrates the case
in which the region $\barr$ is a double polar cap symmetric about the
polar axis $\hat{\mbf z}$. The angular radius of the polar caps is
denoted by $\Theta$. The double polar cap is representative of the
geodetic case in which $\barr$ is the so-called polar gap of missing
observations; its complement $R$ is a latitudinal belt of angular
width $\pi-2\Theta$ around the equator, as shown. In the following,
for brevity, we will shorten all double summations to a notation
requiring only a single sum:
\ber
\sumsh&\rar&\sumshort\hsp,\nnr\\
\sumshL&\rar&\sumshortL\also\sumshR\hspace*{0.43em}
\rar\hspace*{0.7em}\sumshortR.\nnr
\eer

\subsection{Preliminary considerations on the source signal}

We model the geophysical signal as a broadband, square-integrable,
real-valued function $s(\rhat)$ on the surface of the unit sphere
$\Omega=(\theta,\phi)$, defined by the transform pair
\be 
s(\rhat)=\sumshort s_{lm}\Ylm(\rhat),\qquad
s_{lm}=\into s(\rhat)\hsp \Ylm(\rhat)\domg
.
\label{sigdef1}
\ee
The integers $l$ and $m$ are the degree and order of the real
spherical harmonics $\Ylmrh$. These are defined by
\ber
 \Ylm(\theta,\phi)&=&\left\{
\begin{array}{l@{\quad\mbox{if}\hspace{0.6em}}l}
\rule[-2mm]{0mm}{6mm}\sqrt{2}\Xlmth\cos m\phi & -l\le m<0\\
\rule[-2mm]{0mm}{6mm}X_{l0}(\theta)                     & m=0\\
\rule[-2mm]{0mm}{6mm}\sqrt{2}\Xlm(\theta)\sin m\phi & 0< m\le l,\\
\end{array}
\right.\label{Ylm}
\eer
\be
\Xlm(\theta)=
(-1)^m
\sqrt{\frac{\Alma}{4\pi}}\,\Plm
(\cos\theta), 
\label{xlm}
\ee
\be
\Plm (\mu)=
\frac{1}{2^ll!}(1-\mu^2)^{m/2}\left(
\frac{d}{d\mu}\right)^{l+m}\!(\mu^2-1)^l,
\label{plm}
\ee
where $\Plm(\mu)$ is the associated Legendre function, and the
normalization constant \cite[e.g.][]{Edmonds96,Dahlen+98}
\be 
\Alma=\Alm
\label{remind}
.
\ee
With these choices the harmonics $\Ylm(\theta,\phi)$ are
orthonormalized on the unit sphere:
\be
\into\Ylmrh\Ylmprh\domg=\dlmp
.
\label{normalization}
\ee
The fixed-order orthogonality relation for
$\Xlmth$ is
\be
\int_0^{\pi}\Xlm\Xlpm\sin\theta\,d\theta
=\frac{1}{2\pi}\,\dllp.
\label{Xlmortho}
\ee
The addition theorem expresses the sum over all orders of spherical
harmonics at different positions in terms of the angular distance
 $\Delta=\arccos(\rhat\cdot\rhat')$ between them as
\be
\suml_{m=-l}^l\Ylmrh\Ylmrhp =
\tlofp P_l(\rhat\cdot\rhat')
,
\label{additionSH}
\ee
where we note that $P_{l0}=P_0$ and $P_l(0)=1$.  The delta function
$\delta(\rhat,\rhat')=
(\sin\theta)^{-1}\delta(\theta-\theta')\delta(\phi-\phi')$ defined by
\be 
\delta(\rhat,\rhat')=\sum_{l=0}^{\infty}\tlofp\!P_l(\rhat\cdot\rhat')
\label{dirac}
\ee
has the usual sifting property
\be
\into\delta(\rhat,\rhat')\hspace{-0.1em}f(\rhat)\domg=f(\rhat').
\label{replication}
\ee
The sum over all degrees to infinity of the fixed-order colatitudinal
functions at different arguments results in the colatitudinal delta
function: 
\be
(\sin\theta)^{-1}\delta(\theta-\theta')
=2\pi\sum_{l=m}^{\infty}\Xlmth\Xlm(\theta')
.\label{delta}
\ee

\subsection{Noisy measurements at satellite altitude}

For convenience we separate the signal into a bandlimited portion
restricted to the degrees $l=0\rar L$ and a portion over the degrees
$l=L+1\rar\infty$ that complement it:
\be
s(\rhat)=\sumshortL s_{lm}\Ylm(\rhat)+\sumshortR
s_{lm}\Ylm(\rhat)
\label{sigdef2}
,
\ee
where we define $L$ to be the spherical harmonic bandwidth. At the
satellite altitude $a$ above the unit sphere the analytic signal is
given by
\be 
\siup(\rhat)=\sumshortL\suplm\Ylm(\rhat)+\sumshortR\suplm\Ylm(\rhat)
,
\label{sigdef3}
\ee
where the upward continued signal coefficients are given in terms of
the source-level terms by \cite[]{Stacey92,Blakely95}
\be
\suplm=(1+a)^{-l-1} s_{lm}
\label{slmdef3}
.
\ee
The data over the region of coverage $R$ would be given by
eq.~(\ref{sigdef3}) were it not that they are contaminated by
noise. In the uncovered areas $\barr$, no measurements are
available. A satellite thus observes
\ber
d(\rhat)=\left\{\begin{array}{ll}
\siup(\rhat)+n(\rhat) & \mbox{if $\rhat\in R$}\\
\mbox{unknown} & \mbox{if $\rhat\in \bar{R}$}.
\end{array}\right.
\label{datdefs}
\eer
We will restrict attention to the case in which the measurement noise
$n(\rhat)$ is additive and given by a zero-mean stochastic process,
which we assume to be white:   
\ber
\langle n(\rhat)\rangle&=&0,\label{zmean}\\
\langle n(\rhat)n(\rhat')\rangle&=&N\delta(\rhat,\rhat')\label{white}.
\eer
Thus, the power of the noise is denoted by $N$, and we use angular
brackets to denote the ensemble averaging over all possible
realizations required to define the process mean and its spatial
(co)variance.

Combining eqs~(\ref{sigdef3})--(\ref{slmdef3}) and~(\ref{sigdef1})
with the definition~(\ref{additionSH}), we can write the signal
observed at orbital level as a convolution of the surface-level signal
in the form
\be
\siup(\rhat)=\into \Gamma(\rhat,\rhat')s(\rhat')\,\domg'
\label{spread}
,
\ee
where we have defined a ``point spread function''
\be
\Gamma(\rhat,\rhat')=
\sum_{l=0}^{\infty}(1+a)^{-l-1}\tlofp P_l(\rhat\cdot\rhat')
\label{psf}
.
\ee
Thus, the value of the potential field that is observed at a point
$\rhat$ outside the unit sphere is a weighted mixture of the function
values at $\rhat$ and distant other points $\rhat'$ on the unit
sphere. Measurements taken by a satellite at $a>0$ are affected by
regions it does not fly over directly. A satellite thus does probe
into the uncovered regions; conversely, in regions of coverage, it may
be affected by uncovered areas.  As the fractional altitude $a$
increases, the convolution kernel $\Gamma(\rhat,\rhat')$ is
increasingly supported globally. On the other hand, when $a=0$,
eq.~(\ref{psf}) returns the delta function, eq.~(\ref{dirac}), and
eq.~(\ref{spread}) merely illustrates its sifting
property~(\ref{replication}).

\subsection{A new basis for bandlimited field estimators}

We seek an estimate $\hat{s}(\rhat)$ of the signal in
eq.~(\ref{sigdef1}), at the level of the source, from the data
$d(\rhat)$, given by eq.~(\ref{datdefs}), at altitude $a$. It is
crucial to realize that, although any real physical signal $s(\rhat)$
will in general be infinite-band, our estimate $\hat{s}(\rhat)$ must
always be bandlimited. We are thus at liberty to define a new,
bandlimited set of basis functions, to replace the spherical
harmonics. In this manner, the broadband source field can be expressed
as
\be
s(\rhat)=\sumapot s_\alpha g_\alpha(\rhat)+\sumshortR s_{lm}\Ylmrh
.
\label{sigdefx}
\ee
whereas the bandlimited estimated field is given by
\be
\hat{s}(\rhat)=\sumshortL\hat{s}_{lm}\Ylm(\rhat)=\sumapot\hat{s}_\alpha
g_\alpha(\rhat)
.
\label{shatdef1}
\ee
These new bandlimited basis functions $g_\alpha(\rhat)$ indexed by
$\alpha=1,\dots,\Lpot$ will themselves be combinations of spherical
harmonics, inasmuch as they are defined by the transform pair
\begin{subequations}
\label{inversion}
\ber
g_\alpha(\rhat)&=&\sumshortL \galm\Ylm(\rhat),\\
\Ylm(\rhat)&=&\sumapot \glma g_\alpha(\rhat),\\ 
\galm&=&\into g_\alpha(\rhat)\hsp \Ylm(\rhat)\domg
.
\eer
\end{subequations}
The new basis will be rendered orthonormal by requiring that
\begin{subequations}
\label{gortho}
\ber
\into g_{\alpha}(\rhat)g_{\beta}(\rhat)\domg&=&\dab,\\
\sumshortL\galm\glmb&=&\dab,\\
\sumapot\glma\galmp&=&\dllp\dmmp
.
\eer
\end{subequations}
The transformation of the spherical harmonic basis coefficients
$\hat{s}_{lm}$ of the estimate to the expansion coefficients
$\hat{s}_\alpha$ in the new basis is achieved by
\be
\hat{s}_{lm}=\sumapot\glma\hat{s}_{\alpha}\also
\hat{s}_\alpha=\sumshortL\galm\hat{s}_{lm}
\label{rota}
,
\ee
as can be easily deduced by combining eq.~(\ref{shatdef1}) with the
orthonormality conditions, eqs~(\ref{normalization}) and~(\ref{gortho}),
and using eq.~(\ref{inversion}). 

Upward continued to the satellite altitude $a$, the estimate in either
basis is
\be 
\siuph(\rhat)=\sumshortL
\suplmh\Ylm(\rhat)=\sumapot\hsiupa g_\alpha(\rhat)
,
\label{shatdef2}
\ee
where the spherical harmonic coefficients are naturally given by 
\be 
\suplmh=(1+a)^{-l-1} \hat{s}_{lm}
,
\label{sigdef4}
\ee
and the upward continued expansion coefficients of the new basis by 
\begin{subequations}
\label{siupa}
\ber
\hsiupa&=&\sumshortL \galm (1+a)^{-l-1}\hat{s}_{lm}\\
&=&\sumbpot\left(\sumshortL \galm (1+a)^{-l-1}
\glmb\hsp \right)\hat{s}_\beta
,
\eer
\end{subequations}
as is verified by combining eqs~(\ref{rota}) and ~(\ref{sigdef4}). As
expected, eq.~(\ref{siupa}) reduces to a trivial identity at the
surface of the unit sphere, i.e. when $a=0$, by virtue of
eqs~(\ref{gortho}) or~(\ref{rota}).

Given that the measurements made by the satellite are restricted to
the domain $R$ on the unit sphere $\Omega$, we consider it natural to
require of the new basis functions that they be optimally concentrated
on this domain.  We will seek a basis of functions $g(\rhat)$ whose
energy is maximally concentrated inside of the domain $R$ by
maximizing the spatial energy ratio,
\be
\lambda=\fracd{\intr
g^2\domg}{\into^{}g^2\domg}=\mbox{maximum}.
\label{normratio}
\ee
Eq.~(\ref{normratio}) is a statement of Slepian's problem, a classic
in one-dimensional time-series analysis
\cite[]{Slepian83,Percival+93}, on the two-dimensional sphere, 
which we have recently studied in detail \cite[]{Simons+2006a}. In the
next section, we review the main properties of the general solution of
eq.~(\ref{normratio}) for concentration domains of arbitrary geometry,
which we subsequently specialize to the geodetic context by imposing
the circular symmetry of the double-cap polar gap.

\section{S~l~e~p~i~a~n'~s{\hsps}s~p~h~e~r~i~c~a~l{\hsps}p~r~o~b~l~e~m}  
\label{ssp}
\label{secsym0}

In Slepian's problem, the concentration of a bandlimited function $g$
given by 
\be 
g=\sumshortL g_{lm}\Ylm,\qquad
g_{lm}=\into g\hsp \Ylm\domg
,
\label{bandlg}
\ee
to a region $R$ of area $A$ on the unit sphere $\Omega$ is expressed
as the norm ratio, eq.~(\ref{normratio}). Maximization of this
concentration criterion can be achieved in the spectral domain by
solving the algebraic eigenvalue problem
\be
\sfD\hspace{0.05em}\sfg=\lambda\sfg
,
\label{eigen1}
\ee
where $\sfg$ is the $\Lpot$--dimensional spherical harmonic
coefficient column vector 
\be
\sfg=(g_{00}\,\cdots\,g_{lm}\,\cdots\,g_{LL})\T
\label{sfgdef}
\ee
and $\sfD$ is the $\Lpot\times\Lpot$--dimensional matrix 
\be
\sfD = \left(\begin{array}{ccc}
D_{00,00} & \cdots & D_{00,LL} \\
\vdots & {} & \vdots \\
D_{LL,00} & \cdots & D_{LL,LL}
\label{Dmatrixdef}
\end{array}\right),
\ee
whose elements $\Dlmlmp$, $0\le l\le L$ and $-l\le m\le l$, are  
\be
\Dlmlmp=\intr\Ylm\Ylmp\domg
.
\label{Dlmlmpdef}
\ee
The obvious symmetry $\sfD\T=\sfD$ guarantees that the
eigenvectors $\sfg_1,\sfg_2,\ldots,\sfg_{(L+1)^2}$
are mutually orthogonal. We choose them to be orthonormal: 
\be
\sfg_{\alpha}\T\sfg_{\beta}=\dab
\also
\sfg_{\alpha}\T\sfD\hspace{0.05em}\sfg_{\beta}
=\lambda_{\alpha}\dab.
\label{orthogspec}
\ee
The resulting Slepian functions $g$ are orthonormal over the whole
sphere $\Omega$ and orthogonal over the region $R$:
\be
\into g_{\alpha}g_{\beta}\domg=\dab\also
\intr  g_{\alpha}g_{\beta}\domg=\lambda_{\alpha}\dab.
\label{orthogspace}
\ee
The leftmost equations in eqs~(\ref{orthogspec})--(\ref{orthogspace})
correspond to the conditions of eq.~(\ref{gortho}) and guarantee
that the solution indeed forms a valid orthonormal basis. The rightmost
equations illustrate the so-called double orthogonality of the Slepian
basis \cite[]{Gilbert+77}, which, as we will see in a later section,
is a central feature of their utility for the geodetic estimation
problem. 

An approach equivalent to the maximization of eq.~(\ref{normratio}) is
to find broadband functions $h(\rhat)$ that are spacelimited to the
domain $R$, but spectrally concentrated in a bandwidth interval
$0\le l\le L$. The concentration measure in this case,
\be
{\lambda}=
\fracd{\sumshortL h^2_{lm}}
{\sumshort h^2_{lm}}=\mbox{maximum},
\label{normratio2}
\ee
is satisfied by the eigenfunctions of a Fredholm integral eigenvalue
equation in the spatial domain:
\be
\intr \Drhrhp \, h(\rhat')\domg'=
\lambda\hspace{0.05em}h(\rhat),
\quad\quad\rhat\in R
.
\label{firsttimeint}
\ee
The symmetric kernel of eq.~(\ref{firsttimeint}) depends only on
the geodesic angular distance, $\Delta$, between $\rhat$ and $\rhat'$:
\be
\Drhrhp=\sum_{l=0}^L\tlofp\!P_l(\rhat\cdot\rhat') 
.
\label{banddelta}
\ee

The problems of finding bandlimited functions $g$
concentrated to a spatial interval or spacelimited functions $h$ 
concentrated in a spectral interval are completely equivalent. The
domain of eq.~(\ref{firsttimeint}) can be extended to the entire
sphere in which case it applies to the bandlimited functions $g$:
\be
\intr  \Drhrhp \,\grhp\domg'=
\lambda\hspace{0.05em}g(\rhat),
\quad\rhat\in\Omega,
\label{secondtimeint}
\ee
We normalize such that the eigenfunctions $g$ that maximize the
spatial energy ratio~(\ref{normratio}) are identical, within the region
$R$, to the eigenfunctions $h$  maximizing the spectral
ratio~(\ref{normratio2}): 
\be
h(\rhat)=\left\{\begin{array}{ll}
g(\rhat) & \mbox{if $\rhat\in R$}\\
0 & \mbox{otherwise}.
\end{array}\right.
\label{hequalsg}
\ee
The relation
\be
h_{lm}=\sumshortLp \Dlmlmp g_{l'm'}
\label{hlm}
\ee
expresses the coefficients $h_{lm}$, where $0\leq l\leq\infty$, in
terms of the coefficients $g_{lm}$, with $0\leq l\leq L$. This is a
straightforward consequence of the definitions in eqs~(\ref{bandlg}),
(\ref{Dlmlmpdef}) and~(\ref{hequalsg}), and, by eq.~(\ref{eigen1}), it
amounts to $h_{lm}=\lambda\hspace{0.05em}g_{lm}$ when $0\leq l\leq L$.
The eigenvalues of eqs~(\ref{eigen1}) or (\ref{firsttimeint}),
\be 
1>\lambda_1\ge\lambda_2\,\cdots\,\ge\lambda_{\Lpot }>0
,
\label{eigorder}
\ee
measure the quality of the spatiospectral concentration: the
bandlimited function that is most concentrated inside $R$ is $g_{1}$,
with $\lambda_{1}$ being the largest associated eigenvalue, and so
on. The sum of the eigenvalues, or Shannon number, equal to the trace
of $\sfD$, defines a diagnostic area-bandwidth product
\be
K=\sum_{\alpha =1}^{\Lpot}\lambda_{\alpha}=
\intr D(\rhat,\rhat)\domg=\Lpot\,\frac{A}{4\pi}
.
\label{tracedef}
\ee
Spherical Slepian functions of equal Shannon number are scaled
versions of each other in the asymptotic limit $A\rightarrow 0$ and
$L\rightarrow\infty$ with $K$ held fixed
\cite[]{Simons+2006a}. Whenever the area $A$ of the region $R$ is a
small fraction of the area of the sphere, $A\ll 4\pi$, i.e., when
$K\ll\Lpot$, there will be many more well excluded eigenfunctions with
insignificant eigenvalues ($\lambda\approx 0$) than well concentrated
eigenfunctions with significant eigenvalues ($\lambda\approx 1$). If
on the other hand, $R$ covers most of the sphere so that $A\approx
4\pi$ and $K\approx\Lpot$, there will be many more well concentrated
eigenfunctions than well excluded ones.

The sum of the squares of the $\Lpot$ bandlimited eigenfunctions
$g(\rhat)$ is independent of position on the sphere:
\be
\sumapot g_{\alpha}^2(\rhat)=\frac{\Lpot}{4\pi}=\frac{K}{A}.
\label{sumofsq}
\ee
Since the first $K$ eigenfunctions $g_1,g_2,\ldots,g_K$ have eigenvalues
near unity and lie mostly within $R$, and the remainder
$g_{K+1},g_{K+2},\ldots,g_{(L+1)^2}$ have eigenvalues near zero and
lie mostly in the complementary region $\bar{R}=\Omega-R$, the
eigenvalue-weighted sum of squares is well approximated by
\be
\sumapot\lambda_{\alpha}\hspace{0.05em}
g_{\alpha}^2(\rhat)\approx
\left\{\begin{array}{ll}
K/A & \mbox{if $\rhat\in R$}\\
0 & \mbox{otherwise.}
\end{array}\right.
\label{sumofsq2}
\ee
The terms with $K<\alpha\le\Lpot$ should be negligible, so it is
immaterial whether they are included in the sum~(\ref{sumofsq2}) or
not. Taken together, the first $K$ orthogonal eigenfunctions
$g_{\alpha},\alpha=1,2,\ldots,K$, with significant eigenvalues
$\lambda_{\alpha}\approx 1$, provide an essentially uniform coverage
of the \mbox{region $R$}. Rather than requiring $\Lpot$ basis
functions to represent an arbitrary spatially concentrated bandlimited
function, the first $K=\Lpot A/(4\pi)$ members of the
Slepian basis provide a very reasonable approximation.

We shall denote the operator localizing to the complementary region
$\bar{R}$ by $\sfbD$, its eigenfunctions by $\bar{g}$ and its
eigenvalues by $\bar{\lambda}$. It follows from the orthogonality
relation~(\ref{normalization}) that the elements of $\sfbD$ are
\begin{subequations}
\label{blabla}
\ber
\bDlmlmp&=&\intbr\Ylm\Ylmp\domg\label{bDlmlmpdefrepeat}\\
&=&\dlmp-D_{lm,l'm'}
.
\eer
\end{subequations}
The eigenfunctions of $\sfbD$ are identical to those of $\sfD$, but
their ordering indices are reversed. The bandlimited function that is
most concentrated within $\bar{R}$ is most excluded from $R$, i.e.
$\bar{g}_1=g_{\Lpot}$, with an associated eigenvalue
$\bar{\lambda}_1=1-\lambda_{\Lpot}$, and so on.

The localization operator $\sfD$ has an inverse satisfying
\be
\sumshortLpp\Dlmlmpp^{-1}\Dlmpplmp^{}=\dllp\dmmp
,
\label{inverse}
\ee
and for which $\sfD^{-1}\sfg=\lambda^{-1}\sfg$. For future reference,
the inverse of a weighted sum of the localization matrix and its
complement obeys
\begin{subequations}
\label{invfut}
\ber
(\sfD+\eta\sfbD)^{-1}\hspace{0.05em}\sfg
&=&\left[\lambda+\eta(1-\lambda)\right]^{-1}\sfg,\\
\sfg_{\alpha}\T(\sfD+\eta\sfbD)^{-1}\hspace{0.05em}\sfg_{\beta}
&=&\left[\lambda_\alpha+\eta(1-\lambda_\alpha)\right]^{-1}\dab
,
\eer
\end{subequations}
for any weighting parameter $\eta$. Finally, we may extend the
validity of the definition~(\ref{Dlmlmpdef}) to include all degrees
$l\le \infty$ and use eqs~(\ref{additionSH})--(\ref{replication}) to
prove that
\be
\sumshortpp\Dlmlmpp\Dlmpplmp=\Dlmlmp\label{axi0}
.
\ee

\section{A~x~i~s~y~m~m~e~t~r~i~c{\hsps}d~o~m~a~i~n~s}
\label{secpc}

In the previous section we showed that the optimally concentrated
bandlimited basis functions that are the solutions to the Slepian
problem are found by numerical diagonalization of the operator in
eq.~(\ref{Dlmlmpdef}). That this is in general possible for arbitrary
geometries was shown by \cite{Simons+2006a}. However, the particular
geometry of data acquisition on the sphere in the geodetic estimation
problem (Figure~\ref{sddiagram}) allows for substantial
simplifications of this general result. We discuss the special case of
finding concentrated basis functions on the latitudinal belt, the
domain over which satellite measurements are made, via the
concentration within the single and the double polar cap. As we have
seen, the eigenfunctions on a domain $R$ are identical to those on a
complementary spherical domain $\bar{R}$, but with their ordering
indices reversed. Identifying $R$ with the polar caps, rather than
their complement, the belt, as we do -- {\it in this section and the
one that follows only} -- greatly simplifies the equations.

\subsection{Concentration within an axisymmetric polar cap}
\label{secsym1}
When the region of concentration is a circularly symmetric cap of
colatitudinal radius $\Theta$, centered on the north pole, i.e.
\be
R=\big\{\theta: 0\le\theta\le\Theta\big\}
,
\label{polarcap}
\ee
of area $A=2\pi(1-\cos\Theta)$, the matrix elements of eq.~(\ref{Dlmlmpdef}) reduce to
\be
\Dlmlmp =
2\pi\,\dmmp\int_{0}^{\Theta}
\Xlm\Xlpm\sin\theta\,d\theta
.
\label{blockdia}
\ee
The Kronecker delta $\delta_{mm}$ renders the matrix $\sfD$ of
eq.~(\ref{Dmatrixdef}) block-diagonal,
\be
\sfD={\rm diag}\,\left(\sfD^0,\sfD^1,\sfD^1,\ldots,
\sfD^L,\sfD^L\right),
\label{blockdia2}
\ee
where every submatrix $\sfD^m\not=\sfD^0$ occurs twice due to the
doublet degeneracy of $\pm m$. Rather than solving the complete
$\Lpot\times\Lpot$ eigenvalue equation~(\ref{eigen1}), we may solve a
series of $(L-m+1)\times (L-m+1)$ spectral-domain eigenvalue problems,
one for each non-negative order $m$, 
\be
\sfD\sfg=\lambda\sfg
\label{separate},
\ee
where we have dropped the superscript identifying the order. The
eigenvalues belonging to every nonzero order, $m>0$, occur twice.  In
eq.~(\ref{separate}) the column vector $\sfg$ collects the spherical
harmonic coefficients of order $m$:
\be
\sfg=(g_{m}\,\cdots\,g_{l}\,\cdots\,g_{L})\T
,
\label{sfgdef2}
\ee
and the fixed-order matrix $\sfD$ is of the form 
\be
\sfD = \left(\begin{array}{ccc}
D_{mm} & \cdots & D_{mL}\\
\vdots && \vdots\\
D_{Lm} & \cdots & D_{LL}\end{array}\right)
,
\label{littleD&g}
\ee
where, for a particular order $0\leq m\leq L$,
\be
D_{ll'}=2\pi\int_{0}^{\Theta}\Xlm\Xlpm\sin\theta\,d\theta.
\label{kernel4}
\ee
Various methods exist to evaluate the elements of eq.~(\ref{kernel4})
\cite[]{Wieczorek+2005,Simons+2006a}. The important point is that,
while symmetric, and banded, the matrix $\sfD$ is never sparse.  Its
construction thus requires on the order of $(L-m+1)^2/2$ integrals
each.

We rank the $L-m+1$ eigenvalues
$\lambda_1,\lambda_2,\ldots,\lambda_{L-m+1}$ obtained by solving the
fixed-order problem~(\ref{separate}) so that
\be
1>\lambda_1\ge \lambda_2\ge\,\cdots\,\ge\lambda_{L-m+1}>0,
\label{fixedmorder}
\ee
and orthonormalize the eigenvectors $\sfg_1$, $\sfg_2,
\ldots,\sfg_{L-m+1}$ as in eq.~(\ref{orthogspec}). The
associated bandlimited eigenfunctions
$g_1(\theta),g_2(\theta),\ldots,g_{L-m+1}(\theta)$, are given by
\be
g=\sum_{l=m}^Lg_{l}\Xlm,\qquad 
g_l=2\pi\int_0^{\pi}g\hsp \Xlm\sin\theta\,d\theta
,
\label{polarg}
\ee
and satisfy the colatitudinal orthogonality relations
\begin{subequations}
\label{fixedmortho}
\ber
2\pi\int_0^{\pi}g_{\alpha}g_{\beta}\sin\theta\,d\theta&=&\dab,\\ 
2\pi\int_0^{\Theta}g_{\alpha}g_{\beta}\sin\theta\,d\theta 
&=&\lambda_{\alpha}\dab.
\eer
\end{subequations}

The optimally concentrated spatial eigenfunctions $g(\rhat)$
for a given order $-L\leq m\leq L$ are expressed in terms of the
fixed-order colatitudinal eigenfunctions~(\ref{polarg}) by
\be
g(\theta,\phi)=\left\{
\begin{array}{l@{\quad\mbox{if}\hspace{0.6em}}l}
\rule[-2mm]{0mm}{6mm}\sqrt{2}\,g(\theta)\cos m\phi & -L\le m<0\\
\rule[-2mm]{0mm}{6mm}g(\theta)                     & m=0\\
\rule[-2mm]{0mm}{6mm}\sqrt{2}\,g(\theta)\sin m\phi & 0< m\le L.\\
\end{array}
\right.
\label{polarg2}
\ee
The fixed-order Shannon number
\be
K_m=\sum_{\alpha =1}^{L-m+1}\lambda_{\alpha}
,
\label{tracedefaxi}
\ee
again is simply the trace of the fixed-order matrix $\sfD$.

We further note that the complementary fixed-order matrices are given
by 
\be
\bar{D}_{ll'}=\dllp-\bDllp.  
\label{kernel5d}
\ee
The eigenfunctions of the fixed-order matrix $\sfbD$ are identical to
those of $\sfD$ but appear in reverse order, and their eigenvalues sum
to one. The axisymmetric inversion formula analogous to
eq.~(\ref{inverse}) is
\be
\sumelpp \Dllpp^{-1}\Dlpplp^{}=\dllp
,
\label{inverse2}
\ee
and the axisymmetric analog to eq.~(\ref{axi0}) is
\be
\sum_{l''=m}^{\infty}\Dllppm\Dlpplpm=\Dllpm
.
\label{axi1}
\ee

\subsection{Concentration within a double polar cap}
\label{secsym2}
When the region of concentration is a pair of axisymmetric antipodal
caps of colatitudinal radius $\Theta$, i.e., when
\be
R=\big\{\theta: 0\le\theta\le\Theta\big\}\cup
\big\{\theta: \pi-\Theta\le\theta\le\pi\big\}
,
\label{doublecap}
\ee
of area $A=4\pi(1-\cos\Theta)$, the reflection symmetry
\be 
\Xlm(\pi-\theta)=(-1)^{l+m}\Xlmth
\label{reflect}
\ee
checkers the fixed-order matrices $\sfD$ with zeroes, following
\be
D_{ll'}=2\pi\left[1+(-1)^{l+l'}\right]
\int_{0}^{\Theta}\Xlm\Xlpm\sin\theta\,d\theta.  
\label{kernel5}
\ee
Comparison of eqs~(\ref{kernel4}) and (\ref{kernel5}) reveals 
that the eigenfunctions of the double-cap problem can be trivially
obtained from the kernels belonging to the single polar cap.

The spherical Slepian functions resulting from the diagonalization of
the double-cap kernel in eq.~(\ref{kernel5}) are either even or odd
across the equator. Indexing their parity by $p$, as even ($p=e$) or
odd ($p=o$), we modify eq.~(\ref{polarg}) to explicitly skip every other
degree by using a primed summation symbol,
\be
g_p=\sumpp g_{l}\Xlm
\label{grep}
,
\ee
where the lower limit $m_p$ is given by
\be
m_e=m \also m_o=m+1\label{mp}
,
\ee
and the upper limit $L_p$ is
\begin{subequations}
\label{Lpdef}\ber
L_e&=&\left\{\begin{array}{ll}
L & \mbox{if $m$ and $L$ have the same parity}\\
L-1 & \hspace{6.75em}\mbox{opposite}
\end{array}\right.
\\
L_o&=&\left\{\begin{array}{ll}
L-1 & \mbox{if $m$ and $L$ have the same parity}\\
L & \hspace{6.75em}\mbox{opposite}
\end{array}\right.
\eer
\end{subequations}
In this formalism, the coefficients that are required for $g_{e}$ are
$g_{m}$, $g_{m+2},\dots$, $g_{L}$ if $m$ and $L$ are both even or both
odd, and $g_{m}$, $g_{m+2},\dots$, $g_{L\!-1}$ if $m$ and $L$ are of
opposite parity. Likewise, the coefficients of $g_o$ are $g_{m+1}$,
$g_{m+3},\dots$, $g_{L-1}$ if $m$ and $L$ are both even or both odd,
and $g_{m+1}$, $g_{m+3},\dots$, $g_{L}$ if $m$ and $L$ have opposite
parity. Eqs~(\ref{reflect}) and (\ref{grep}) then confirm that
\be
g_e(\theta)=g_e(\pi-\theta)\also
g_o(\theta)=-g_o(\pi-\theta)
.
\ee
While an equation of the form~(\ref{separate}) returns an alternation
of even and odd functions with decreasing eigenvalues $\lambda$, the
indices of the matrix $\sfD$ may be permuted to form a block-diagonal
form
\be
\sfD={\rm diag}\,\left(\sfD_e,\sfD_o\right)
,
\ee
for which the half-size separate eigenvalue equations 
\be
\sfD_e\sfg_e=\lambda_e\sfg_e \also 
\sfD_o\sfg_o=\lambda_o\sfg_o
\label{separate3}
\ee
return exclusively even or odd solutions. This avoids round-off
problems and speeds up the diagonalization. 

The slight perversity of our notation is that, in an all-even or
all-odd approach as in eq.~(\ref{separate3}), writing the fixed-order
coefficient $g_{l}$ or indeed any expression involving the spherical
harmonic degree~$l$, always has to be accompanied by the set of
allowable degrees $l$, since $l= m_p, m_p+2, \dots , L_p$, depending
on the parity.

\section{T~h~e{\hsps}m~a~g~i~c{\hsps}o~f{\hsps}c~o~m~m~u~t~a~t~i~o~n}
\label{tmoc}

While conceptually simple, the formalism presented in the previous
section suffers from two important difficulties. First, assembling the
matrices of eqs~(\ref{kernel4}) and (\ref{kernel5}) requires the
calculation of $\mathcal{O}(L^2)$ matrix elements, by numerical
integration or other means
\cite[]{Wieczorek+2005,Simons+2006a}. Second, and more importantly,
when a large number of near-zero eigenvalues is present, e.g. when
$\Theta\rar0$ and the complementary solutions are sought on $\bar{R}$,
the diagonalization is rarely stable, as discussed by
\cite{Albertella+99}. In principle, any orthogonal set of
solutions might suffice to solve the problem at hand, but those
solutions will vary depending on the method of computation. The method
outlined below always produces stable, unique, solutions, and it does
so at a speed which requires only $\mathcal{O}(L)$ algebraic
evaluations to construct the kernels.

\subsection{A commuting operator for the single polar cap}

In the case of the single symmetric polar cap, eq.~(\ref{firsttimeint})
can be rewritten as a series of fixed-order integral equations
\be
\int_{0}^{\Theta}D(\theta,\theta')\,h(\theta')\sin\theta'\,d\theta'
={\lambda}\hspace{0.05em}h(\theta),\quad 0\le\theta\le\Theta
.
\label{Fredholm}
\ee
each with an $m$-dependent, separable, symmetric kernel 
\be
D(\theta,\theta')=2\pi\suml_{l=m}^{L}\Xlmth\Xlm(\theta'). 
\label{dlx}
\ee
Building on the results derived by \cite{Gilbert+77},
\cite{Grunbaum+82} found a second-order differential operator that
commutes with the convolutional integral operator of
eq.\ (\ref{Fredholm}). For any $0\leq m\leq L$, it is of the form  
\be
\ssT=(\cos\Theta-\cos\theta)\nabla_m^2+\sin\theta\frac{d}{d\theta}
-L(L+2)\cos\theta
,
\label{grunbaumop}
\ee
where $\nabla_{\!m}^2=d^2/d\theta^2+\cot\theta\,(d/d\theta)
-m^2(\sin\theta)^{-2}$ is the fixed-order Laplace-Beltrami
operator. The proof of the commutation relation is sketched in
\cite{Simons+2006a}. Since commuting operators have identical
eigenfunctions, the spacelimited, 
fixed-order eigenfunctions $h(\theta)$ can be found by
solving the differential eigenvalue equation
\be
\ssT h(\theta)=\chi\hspace{0.05em}h(\theta), \quad 0\le\theta\le\Theta
\label{grunval},
\ee
where $\chi\not= \lambda$ is the associated Gr\"{u}nbaum
eigenvalue. 


Gr\"{u}nbaum's operator is a Sturm-Liouville operator
\cite[]{Simons+2006a}. Thus, eq.~(\ref{grunval}) has a simple and
easily sorted spectrum, with an infinite number of distinct
eigenvalues $\chi_1<\chi_2<\ldots$ having an accumulation point at
infinity.  The rank orderings of the eigenvalues
$\chi_1,\chi_2,\ldots$ and the spatiospectral concentration factors
$\lambda_1,\lambda_2,\ldots,\lambda_{L-m+1}$ are reversed, so that the
eigenfunction $h_1(\theta)$ associated with the numerically smallest
eigenvalue $\chi_1$, which has no nodes in the polar cap
$0\le\theta\le\Theta$, is the best concentrated fixed-order
eigenfunction; $h_2(\theta)$, which has exactly one node, is the next
best concentrated, and so on.

Extending the domain of eq.~(\ref{grunval}) to the entire domain
$0\le\theta\le\pi$ transforms the unknown functions from the
spacelimited functions $h$ again into the bandlimited functions $g$. 
Eq.~(\ref{grunval}) is then equivalent to the algebraic eigenvalue
equation  
\be
\sfT\hsp\sfg=\chi\sfg,
\label{grunmatrix}
\ee
where $\sfT$ is the $(L-m+1)\times (L-m+1)$ matrix with coefficients
\be
T_{ll'}=2\pi\int_{0}^{\pi}\Xlm(\ssT\hspace{-0.05em}
\Xlpm)\sin\theta\,d\theta
.
\label{grunmatrix3}
\ee
Eqs~(\ref{grunmatrix})--(\ref{grunmatrix3}) are completely equivalent
to eqs~(\ref{separate}) and (\ref{kernel4}). Both matrices $\sfD$ and
$\sfT$ are symmetric, $\sfD=\sfD\T$ and
$\sfT=\sfT\T$. In addition, they commute,
$\sfD\hspace{0.05em}\sfT=\sfT\hspace{0.05em}\sfD$, so they have
identical eigenvectors.  In index notation,
\be
\sum_{n=m}^LD_{ln}T_{nl'}=\sum_{n=m}^LT_{ln}D_{nl'},
\label{G&Dcommute2}
\ee
which can be used as a numerical check.

There are a number of ways to evaluate the elements of the
Gr\"{u}nbaum matrix in eq.~(\ref{grunmatrix3}), but the important
result is that  $\sfT$ is tridiagonal \cite[][]{Simons+2006a}: 
\begin{subequations}
\ber
T_{ll}&=&-l(l+1)\cos{\Theta},\\
T_{l\,l+1}&=&\big[l(l+2)-L(L+2)\big]\\
&&{}\times\sqrt{\fracd{(l+1)^2-m^2}{(2l+1)(2l+3)}},\nnr\\
T_{ll'}&=&0\quad\mbox{otherwise}
.
\eer
\end{subequations}

Eq.~(\ref{grunmatrix}) can be used to find the $(L-m+1)$--dimensional
eigenvectors $\sfg$ and thus the optimally concentrated polar cap
eigenfunctions $g(\theta)$ by numerical diagonalization of a
tridiagonal matrix $\sfT$ with analytically prescribed elements and a
spectrum of eigenvalues $\chi$ that is guaranteed to be
regular. Unlike the diagonalization of the original matrix $\sfD$ in
eq.(\ref{separate}), this procedure enables the stable computation of
bandlimited functions that are optimally concentrated in a large
rather than a small region of the unit sphere, as may be the case in
geodesy.

\subsection{A commuting operator for the double polar cap}

Knowing that the solutions to the concentration problem for the double
polar cap are either even or odd across the equator, we may write the
integral equation~(\ref{firsttimeint}), by analogy with
eq.~(\ref{Fredholm}), as follows. Indicating the parity of the
solutions by the subscript $p$, which takes the values $p=e$ for the
even solutions and $p=o$ for the odd solutions, it can be seen that
\be
\int_{0}^\Theta D_p(\theta,\theta')\,h_p(\theta')\sin\theta'\,d\theta'=
\lambda\hspace{0.05em}h_p(\theta),
\label{Fredholm2}
\ee
which is valid inside the double antipodal polar cap
\be
\label{thvalid}
\big\{\theta: 0\le\theta\le\Theta\big\}\cup
\big\{\theta: \pi-\Theta\le\theta\le\pi\big\},
\ee
and where the $m$-dependent kernel, analogous to eq.~(\ref{dlx}), is
\be
D_p(\theta,\theta')=4\pi \sumpp  \Xlmth\Xlm(\theta')
.
\label{Dproof1}
\ee
As in eq.~(\ref{grep}), the primed summation skips every second
entry, and the lower and upper limits are as in
eqs~(\ref{mp})--(\ref{Lpdef}).    

Again basing ourselves on the results of \cite{Gilbert+77} and
\cite{Grunbaum+82}, we show in Appendix~\ref{comm} that a
Sturm-Liouville second-order differential operator that commutes with
the convolutional integral operators of eq.~(\ref{Fredholm2}) is of
the form
\ber
\ssT_p&=&(\cos^2\Theta-\cos^2\theta)\nabla_m^2+
2\cos\theta\sin\theta\frac{d}{d\theta}\nnr\\ 
&&{}-L_p(L_p+3)\cos^2\theta
.
\label{grunbaumop2}
\eer
The individual matrix operators $\sfT_p$ are once again tridiagonal
and symmetric and commute with the even or odd $\sfD_p$ of
eq.~(\ref{separate3}). The bandwidth $L_p$ is the same as
in~(\ref{Lpdef}). The elements of the double-cap Gr\"unbaum matrices are
\begin{subequations}
\label{gdefi}
\ber
T^p_{ll}&=&-l(l+1)\cos^2{\Theta}
+\frac{2}{2l+3}\left[(l+1)^2-m^2\right]\nnr\\
&&{}+[(l-2)(l+1)-L_p(L_p+3)]\nnr\\
&&{}\times\left[\frac{1}{3}-\frac{2}{3}\,
\fracd{3m^2-l(l+1)}{(2l+3)(2l-1)}\right],\\ 
T^p_{l\,l+2}&=&\fracd{\big[l(l+3)-L_p(L_p+3)\big]}{2l+3}\\
&&{}\times\sqrt{\fracd{\left[(l+2)^2-m^2\right]
\left[(l+1)^2-m^2\right]}{(2l+5)(2l+1)}},\nnr\\
T^p_{ll'}&=&0\quad\mbox{otherwise}
.
\eer\end{subequations}
We again emphasize that, since we focus our attention separately on
the kernels returning even or odd eigenfunctions $g_e$ or $g_o$, the
degrees involved are restricted to $l,l'= m_p, m_p+2, \dots ,
L_p$. Since every other degree in the matrix described by
eq.~(\ref{gdefi}) is skipped, both $\sfT_e$ and $\sfT_o$ are
tridiagonal as in the single-cap case. As in eq.~(\ref{separate3}) we
compute the even and odd eigenfunctions by separately solving
\be
\sfT_e\hspace{0.05em}\sfg_e=\chi_e\sfg_e \also 
\sfT_o\hspace{0.05em}\sfg_o=\chi_o\sfg_o
.\label{separate4}
\ee
Subsequently, we establish a single rank order of decreasing
spatiospectral concentration: either per order, as in
eq.~(\ref{fixedmorder}), or across all orders, as in
eq.~(\ref{eigorder}). 

\begin{figure*}
\rotatebox{-90}{
\iftwocol
{\includegraphics[height=0.815\textwidth]{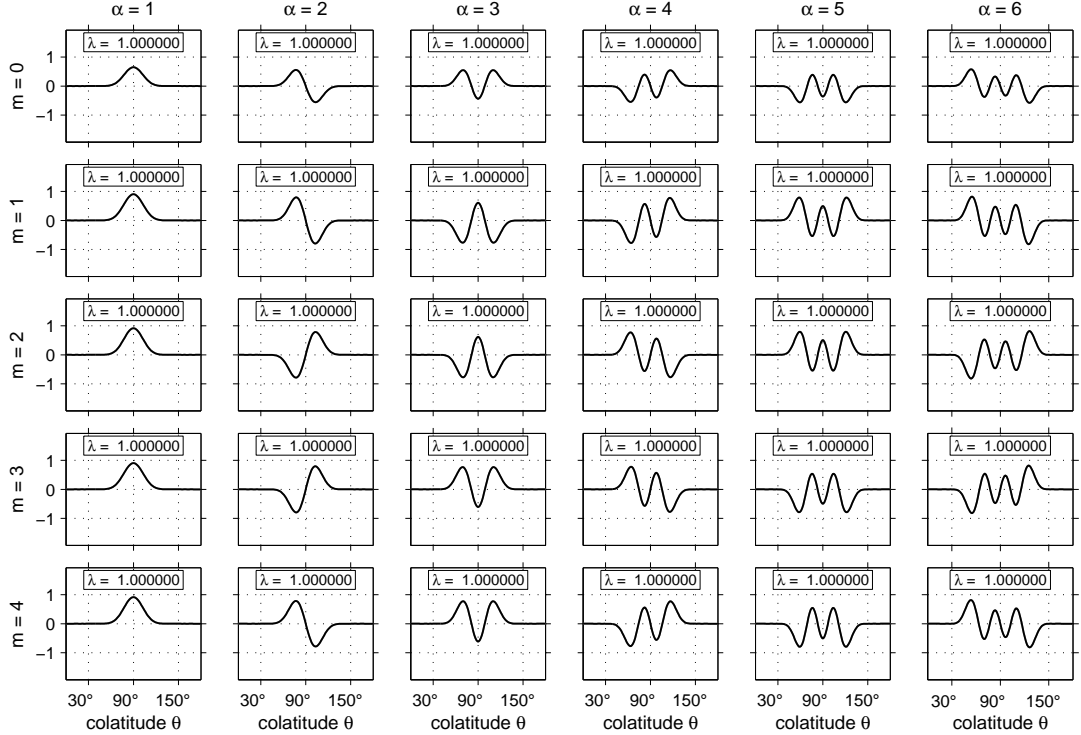}}
{\includegraphics[height=0.85\textwidth]{sdbelt.eps}}
} 
\caption{Colatitudinal dependence of the first six fixed-order,
$m=0\rar4$, eigenfunctions $g_\alpha(\theta)$, $\alpha=1\rar6$,
bandlimited to $L=18$, that are well concentrated in the latitudinal
belt extending $\pm60^{\circ}$ on either side of the equator. The
quality of the spatial concentration is expressed by the labeled
eigenvalues $\lambda_\alpha$. None of the plotted functions show
appreciable energy inside the complementary pair of antipodal polar
caps of radius $\Theta=30^{\circ}$.}
\label{sdbelt}  
\end{figure*} 
\begin{figure*}
\rotatebox{-90}{
\iftwocol
{\includegraphics[height=0.815\textwidth]{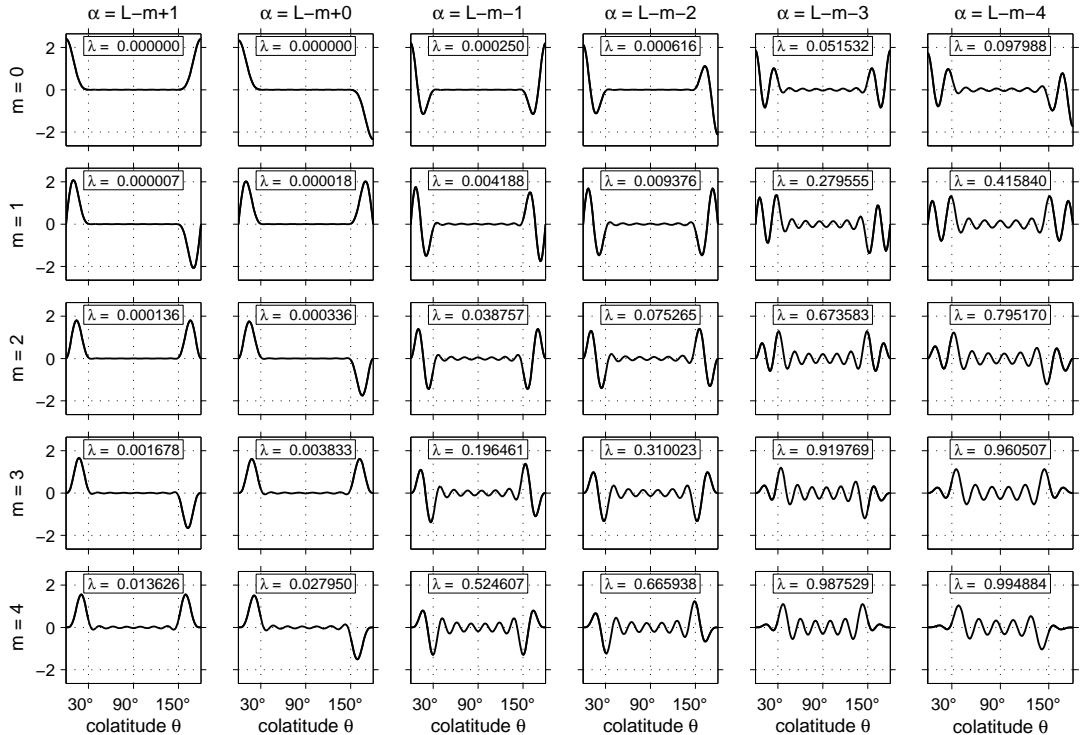}}
{\includegraphics[height=0.85\textwidth]{sdspace.eps}}
} 
\caption{Colatitudinal dependence of the last six fixed-order,
$m=0\rar4$, eigenfunctions $g_\alpha(\theta)$, $\alpha=L-m+1\rar
L-m-4$, bandlimited to $L=18$. These are generally poorly concentrated
in the latitudinal belt $\pm60^{\circ}$ about the equator, except
where the rank $\alpha$ exceeds the fixed-order Shannon number $K_m$
(examples in lower right). The functions that have the least energy
inside of the equatorial belt, as shown by their low eigenvalues
$\lambda_\alpha$, are best concentrated inside the complementary polar
caps of colatitudinal radius $\Theta=30^{\circ}$.  }
\label{sdspace}  
\end{figure*}
\begin{figure*}
\rotatebox{0}{
\iftwocol
{\hspace{-1cm}\includegraphics[width=0.75\textwidth]{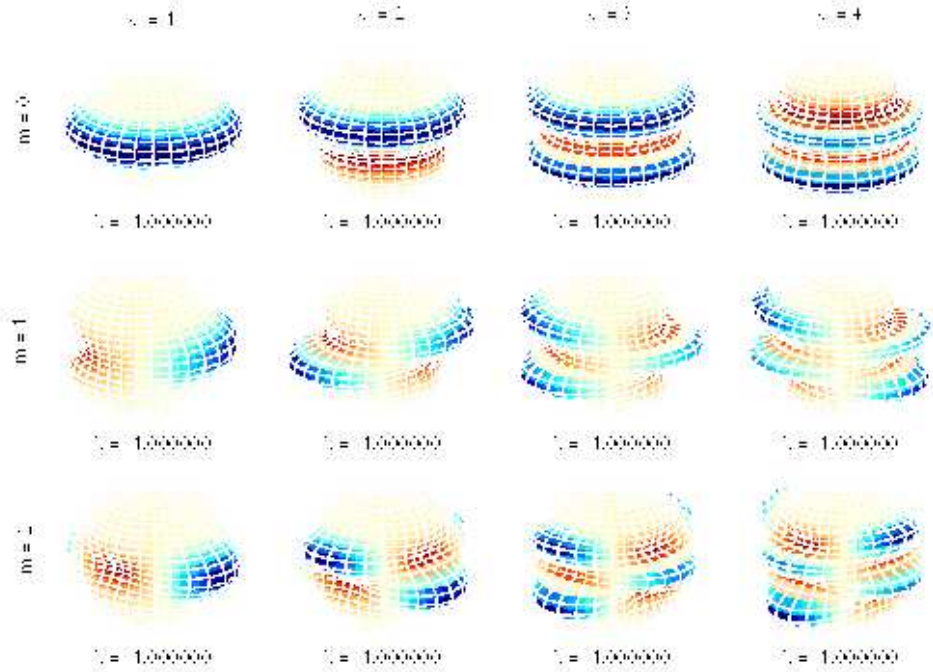}}
{\hspace{-1cm}\includegraphics[width=0.85\textwidth]{sdcmb3.eps}}
} 
\caption{Three-dimensional spatial dependence of the first four fixed-order,
$m=0\rar2$, eigenfunctions $g_\alpha(\theta)$, $\alpha=1\rar4$,
bandlimited to $L=18$, well concentrated in the latitudinal belt
extending $\pm60^{\circ}$ on either side of the equator, as expressed
by their eigenvalues $\lambda_\alpha$. Plot arrangement is as in
Figure\ \ref{sdbelt}.}
\label{sdcmb}  
\end{figure*} 
\begin{figure*}
\rotatebox{0}{
\iftwocol
{\hspace{-1cm}\includegraphics[width=0.75\textwidth]{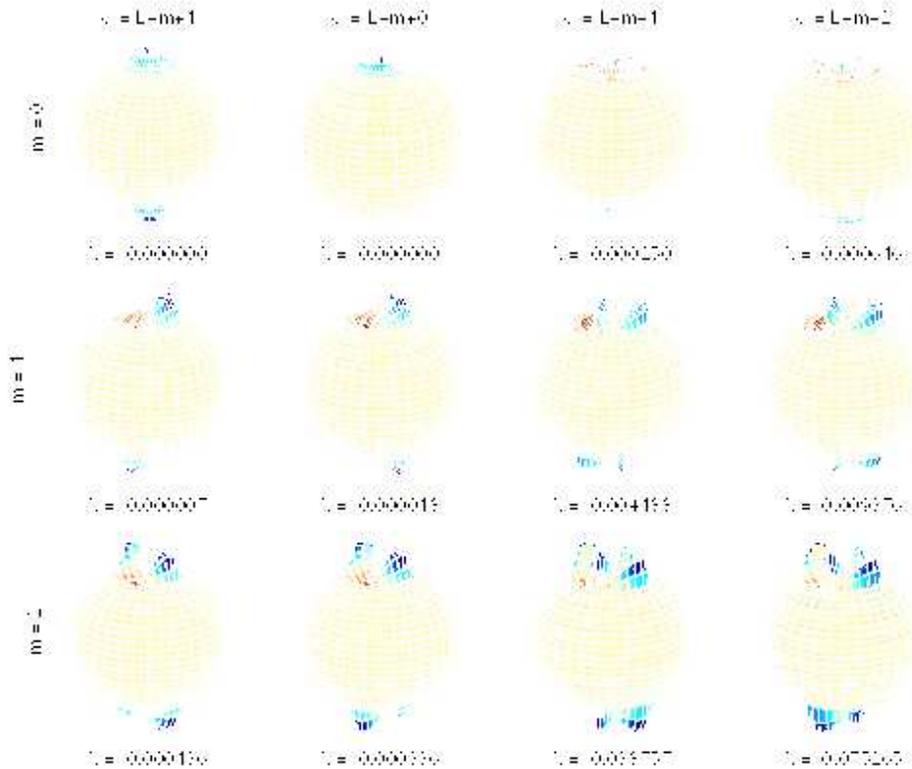}}
{\hspace{-1cm}\includegraphics[width=0.85\textwidth]{sdcmb4.eps}}
} 
\caption{Three-dimensional spatial dependence of the last four fixed-order,
$m=0\rar2$, eigenfunctions $g_\alpha(\theta)$,  $\alpha=L-m+1\rar
L-m-2$, bandlimited to $L=18$, poorly concentrated in the belt
$\pm60^{\circ}$ about the equator, as expressed by their eigenvalues
$\lambda_\alpha$. Plot arrangement is as in Figure\ \ref{sdspace}.}
\label{sdcmb2}
\end{figure*} 

\section{A{\hsps}S~l~e~p~i~a~n{\hsps}b~a~s~i~s{\hsps}o~n{\hsps}t~h~e{\hsps}b~e~l~t}  
\label{eotb}

Concentration within a single polar cap was treated extensively by
\cite{Wieczorek+2005} and \cite{Simons+2006a}. We refer to their
figures for illustrations. In this section we illustrate the solutions
to the concentration problem when the concentration region contains
all but an antipodal pair of polar caps of radius $\Theta$. Reverting
to our notational convention in Section~\ref{ps}, we again use $R$ to
denote an equatorial strip or latitudinal belt extending
$\pi/2-\Theta$ north and south of the equator. Consequently, the
antipodal pair of polar caps themselves is again defined to be the
excluded region $\bar{R}$, in line with their role as the geodetic
polar gap in which no satellite observations are available.

\subsection{Spatial-domain solutions}

The six eigenfunctions $g_\alpha(\theta)$, $\alpha=1\rar6$ of fixed
order $0\leq m\leq 4$ that are most optimally concentrated in the
latitudinal belt complementing a $\Theta=30^{\circ}$ double polar cap
are plotted in Figure~\ref{sdbelt}.  Their associated eigenvalues
$\lambda_\alpha$ are listed to six-figure accuracy. The latitudinal
belt ranges from $60^{\circ}$ north to $60^{\circ}$ south
symmetrically about the equator. With the chosen bandwidth $L=18$, the
Shannon number defined in eq.~(\ref{tracedef}) is
$K=\Lpot\cos\Theta\approx 313$, which approximates the number of well
concentrated eigenfunctions with $\lambda\approx 1$.  The best
concentrated eigensolution of every order is a bell-shaped even
function with no nodes in the belt. In keeping with the
Sturm-Liouville character of the Gr\"unbaum operator, every subsequent
solution acquires one more node, so that the second best of every
order is an odd function, the third is even, and so on. All of the
eigenvalues shown in Figure~\ref{sdbelt}, calculated by numerically
integrating eq.~(\ref{normratio}), are equal to one within six-figure
accuracy, indicating that the concentration to the belt is nearly
perfect, while both poles are almost completely excluded. Since the
concentration region is very large, the calculation of these functions
by any means other than the Gr\"unbaum procedure described above will
fail.

\begin{figure}\center
\rotatebox{0}{
\iftwocol
{\includegraphics[width=1\columnwidth]{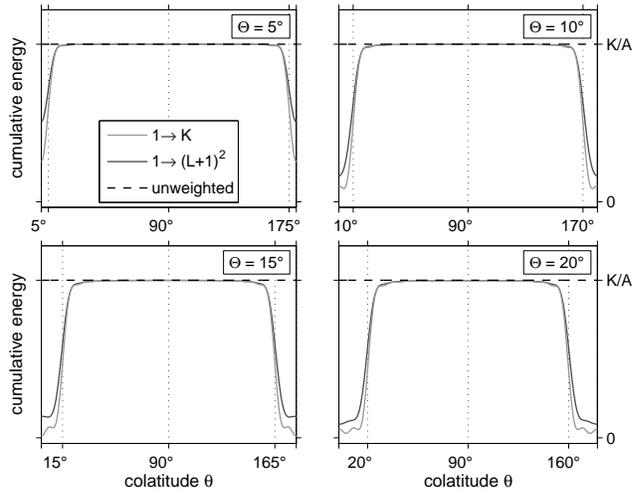}}
{\includegraphics[width=0.85\columnwidth]{sdsumbelt.eps}}
} 
\caption{Cumulative energy of $L=18$ bandlimited eigenfunctions
concentrated inside of belts complementary to antipodal polar caps of
radius $\Theta=5^{\circ},10^{\circ},15^{\circ}$ and $20^{\circ}$. The
Shannon numbers are $K=360,356,349$ and $339$. The  
sums of squares $g_1^2(\theta,\phi)+g_2^2(\theta,\phi)+
\,\cdots$ and
$\lambda_1g_1^2(\theta,\phi)+\lambda_2g_2^2(\theta,\phi)+\,\cdots$ are
plotted versus colatitude $\theta$, along a fixed arbitrary meridian
$\phi$. Dashed lines show the full unweighted sums of $\Lpot$ terms,
which attain the constant value $K/A$ over the entire sphere
$0^{\circ}\leq\theta\leq180^{\circ}$. Solid lines show the
eigenvalue-weighted partial sums of $K$ terms and the full sums of
$\Lpot$ terms, which are very nearly equal, and concentrated uniformly
inside of the belt $\Theta\le\theta\le 180^{\circ}-\Theta$.}
\label{sdsumall} 
\end{figure}
\begin{figure}\center
\rotatebox{0}{
\iftwocol
{\includegraphics[width=1\columnwidth]{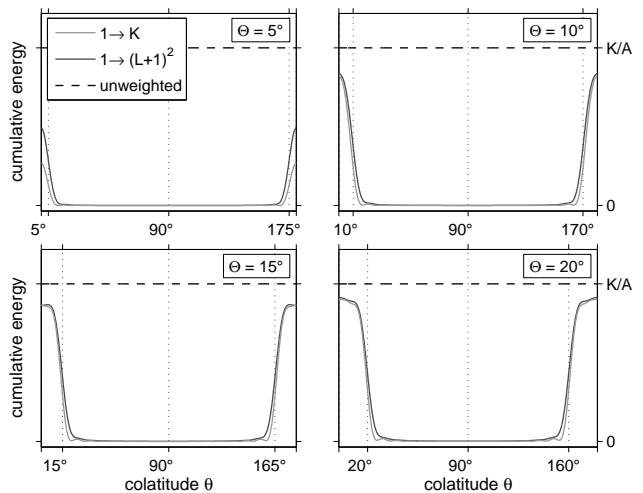}}
{\includegraphics[width=0.85\columnwidth]{sdsumcaps.eps}}
} 
\caption{Cumulative energy of  $L=18$ bandlimited eigenfunctions
concentrated within circularly symmetric polar caps of colatitudinal
radius $\Theta=5^{\circ},10^{\circ},15^{\circ}$ and $20^{\circ}$. The
Shannon numbers are $K=1,5,12,22$. The symbols used are identical to
those of Figure~\ref{sdsumall}. The solid lines showing the
eigenvalue-weighted partial sums of $K$ terms and the full sums of
$\Lpot$ terms are very nearly equal, and concentrated uniformly within
the pair of antipodal caps $0^{\circ}\leq\theta\leq\Theta$ and
$180^{\circ}\hspace{-0.35em}-\Theta\leq\theta\leq180^{\circ}$.}
\label{sdsumall2}
\end{figure}
\begin{figure*}\center
\rotatebox{0}{
\iftwocol
{\includegraphics[width=0.7\textwidth]{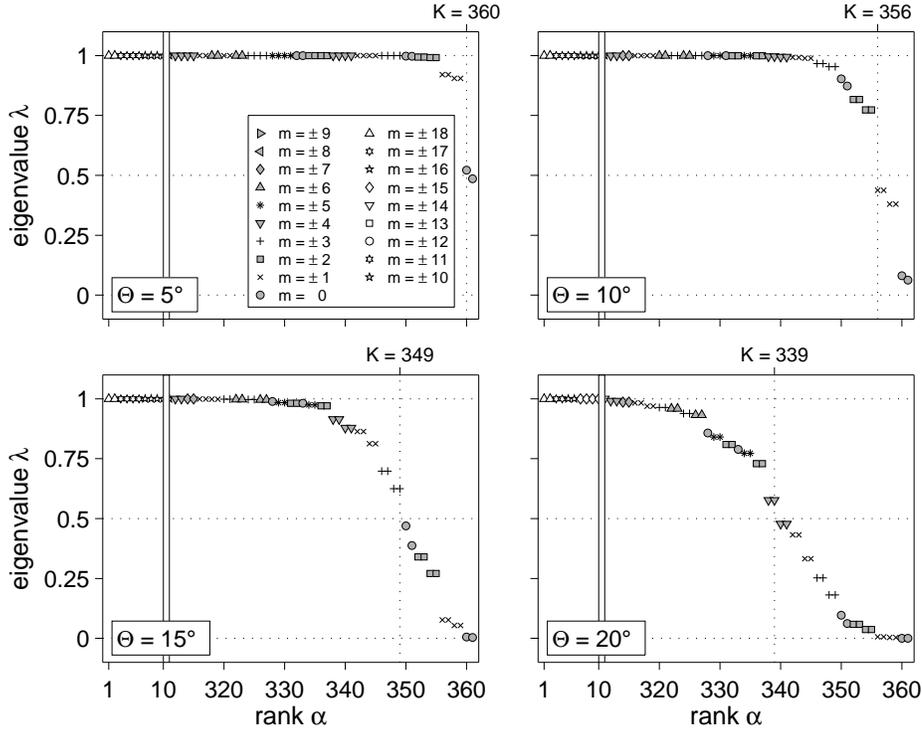}}
{\includegraphics[width=0.85\textwidth]{sdeigs.eps}}
} 
\caption{Eigenvalue spectra of the mixed-order $L=18$
bandlimited operators concentrating within the equatorial belts that
complement antipodal pairs of axisymmetric caps of radius
$\Theta=5^{\circ},10^{\circ},15^{\circ}$ and $20^{\circ}$. The total
number of eigenvalues is $\Lpot=361$; only $\lambda_1\rar\lambda_{10}$
and $\lambda_{312}\rar\lambda_{361}$ are shown. Different symbols are
used for the various orders $-18\leq m\leq 18$; juxtaposed identical
symbols are $\pm m$ doublets. Top labels specify the rounded Shannon
numbers $K=360, 356, 349$ and $339$.}
\label{sdvals}
\end{figure*}

With the parameters unchanged from Figure~\ref{sdbelt},
Figure~\ref{sdspace} shows the six worst concentrated eigenfunctions
on the belt, $g_\alpha(\theta)$, $\alpha=L-m+1\rar L-m-4$. These now
naturally have almost all of their energy inside of the antipodal pair
of polar caps of radius $\Theta=30^{\circ}$. For the zonal functions
of order $m=0$, the even-odd alternation starting at $\alpha=1$ with
an even function in Figure~\ref{sdbelt} ends at $\alpha=L+1$ with an
even function, since $L$ itself is even. At $m=1$, the sequence starts
with an even function but ends at $\alpha=L$ with an odd function, at
$m=2$ with an even function at $\alpha=L-1$, and so on. Thus, the
worst concentrated $m=0$ eigenfunction is even about the equator, the
worst $m=1$ eigenfunction is odd, and so on, in a pattern that
alternates with increasing order. Had $L$ itself been odd, the worst
concentrated zonal function would have been odd, the worst $m=1$
function even, and so on, reversing the pattern.

Three-dimensional perspective views of the first four of the
fixed-order $m=0\rar 2$ functions whose colatitudinal dependence we
plotted in Figure~\ref{sdbelt} are shown in Figure~\ref{sdcmb}. In
accordance with eq.~(\ref{polarg2}) the zonal $m=0$ eigenfunctions do
not display any longitudinal zero crossings, since the number of
longitudinal nodes follows the order $m$. Similarly, in
Figure~\ref{sdcmb2} we plot a three-dimensional rendering of twelve of
the worst concentrated eigenfunctions of Figure~\ref{sdspace}.

In Figure~\ref{sdsumall} we show the eigenvalue-weighted pointwise
sums of squares
$\sum_{\alpha}\lambda_{\alpha}g_{\alpha}^2(\theta,\phi)$ for
latitudinal belts complementary to double polar caps of radii
$\Theta=5^{\circ},10^{\circ},15^{\circ},20^{\circ}$, with a bandwidth
$L=18$. The cumulative sums are concentrated inside of the latitudinal
belt; solid lines in grey and black distinguish the sums carried up to
the first $K$ (the Shannon number) or all $\Lpot$ possible terms. In
contrast, the cumulative sums of the cap eigenfunctions, shown in
Figure~\ref{sdsumall2}, are concentrated within the double polar
cap. The full unweighted sums $\sum_{\alpha}g_{\alpha}^2(\theta,\phi)$
of all $\Lpot$ terms (dashed black lines) are exactly
$K/A=\Lpot/(4\pi)$ over the entire sphere in accordance with
eq.~(\ref{sumofsq}), and the expectation in eq.~(\ref{sumofsq2}) is
confirmed: inside of the concentration domain, the weighted sums
approach $K/A$.

\subsection{Eigenvalue spectra}

In Figure~\ref{sdvals} we show the reordered, mixed-order eigenvalue
spectra for the concentration problem within the latitudinal belt
between polar caps of colatitudinal radii
$\Theta=5^{\circ},10^{\circ},15^{\circ},20^{\circ}$. Once again the
maximal spherical harmonic degree is $L=18$. The rounded Shannon
numbers $K=360,356,349,339$ lie in the middle of the steep,
transitional part of the spectra, roughly separating the reasonably
well concentrated eigensolutions ($\lambda\ge0.5$) from the more
poorly concentrated ones ($\lambda<0.5$) in all four cases. There are
many more functions that are well concentrated in the equatorial strip
than there are that are concentrated inside of the double polar cap,
as shown by the break at $\alpha=10$ in the abscissas.

\begin{figure}\center
\rotatebox{0}{
\iftwocol
{\includegraphics[width=\columnwidth]{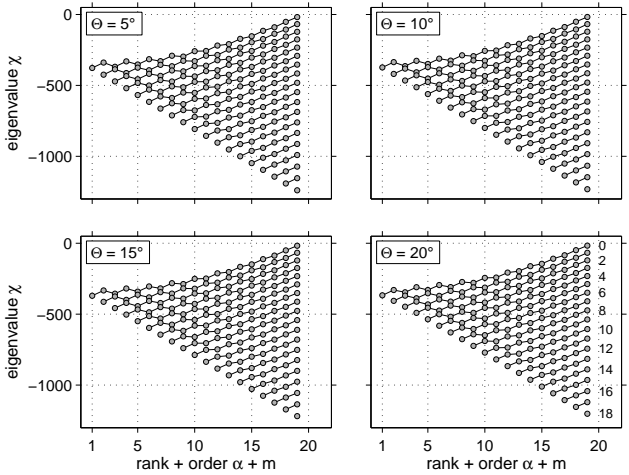}}
{\includegraphics[width=0.85\columnwidth]{sdgrunval.eps}}
} 
\caption{Eigenvalue spectra of the fixed-order $L=18$
bandlimited Gr\"{u}n\-baum operators commuting with the operators whose
eigenvalues are shown in  Figure~\ref{sdvals}. Separate sequences of
eigenvalues $\chi_1,\chi_2, \ldots,\chi_{L-m+1}$ for each angular
order $0\le m\le L$ are connected by lines. Each sequence
is offset horizontally by its order $m$, and vertically by 50 units
per order.}
\label{sdgrunval}  
\end{figure}

The corresponding Gr\"{u}nbaum eigenvalue spectra are shown in
Figure~\ref{sdgrunval}. The ranked eigenvalues $\chi$ for every order $0\leq
m\leq L$ are connected by lines, with each sequence offset
horizontally by its order, and vertically by an arbitrary 50 units, to
facilitate inspection.  Thus, $L+1$ eigenvalues
$\chi_1,\chi_2,\ldots,\chi_{L+1}$ are plotted for $m=0$, whereas a
single eigenvalue $\chi_1$ is plotted for $m=L$. The spacing between
adjacent fixed-order eigenvalues is roughly equant, without the
numerically troublesome plateaus of nearly equal values apparent in
Figure~\ref{sdvals}. This regularity is guaranteed by the
Sturm-Liouville character of the Gr\"unbaum operator $\ssT_p$ in
eq.~(\ref{grunbaumop2}). 

\subsection{Analytic continuation} 
\label{udc}

The Slepian functions $g_{1}(\rhat),\,\dots\,,g_{\Lpot}(\rhat)$ are
defined on the surface of the unit sphere $\Omega$. Together, they
form a natural basis set for the expansion of potential fields and,
in particular, estimates of these fields, on a sphere of radius
$\|\rhat\|=1$, as in eqs.~(\ref{sigdefx})--(\ref{shatdef1}). This new
basis is localized: the support of the first Shannon number $K$ basis
functions lies mostly in the concentration region $R$, whereas the
remainder are concentrated outside of this area of interest, in
$\barr$. With satellite observations we are of course mostly
interested in the signal at some height above the surface of the unit
sphere. We have previously derived an expression for the expansion of
a field estimate at satellite altitude $a$, in eq.~(\ref{siupa}). It
is immediately obvious from this equation that, even if we were only
interested in the first $K$ upward continued Slepian expansion
coefficients of the estimate, we would still need to know and
calculate the full set of $\Lpot$ coefficients at zero altitude.  The
full impact of this statement will not become clear until later in
this paper, but to anticipate it we derive here a set of Slepian basis
functions that are designed specifically to represent signals at an
altitude. We can do this by interpreting the Slepian functions we have
just constructed as potential functions themselves. In that case their
upward harmonic continuation onto a sphere larger radius $\|{\mbf
r}\|=1+a$, where $a>0$, yields new functions $\gup$ for which
\be 
\gup=\sumshortL\guplm\Ylm,\qquad \guplm=(1+a)^{-l-1} g_{lm}
.
\label{bandlgup}
\ee
Had we instead defined the Slepian functions on the larger sphere
to begin with, their analogues downward continued onto the unit sphere
would be obtained as
\be 
\gdn=\sumshortL \gdnlm\Ylm,\qquad \gdnlm=(1+a)^{l+1}g_{lm}
.
\label{bandlgdn}
\ee
It is thus useful to define a symmetric, $\Lpot\times\Lpot$, downward
continuation matrix $\sfA$, whose elements are
\be
\Almlmp=(1+a)^{l+1}\dlmp
,
\label{Adef}
\ee
which allows us to restate the equations relating the upward and downward
continued coefficients to each other concisely as:
\begin{subequations}
\ber
\label{concise}
\sfgup=\sfA^{-1}\sfg,&\qquad& \sfg=\sfA\hsp\sfgup,\\
\sfgdn=\sfA\hsp\sfg,\hspace{0.95em}&\qquad&\sfg=\sfA^{-1}\sfgdn.
\eer
\end{subequations}
The orthogonality relations of eqs~(\ref{gortho})
and~(\ref{orthogspec}) can be rewritten in terms of $\gup$ and $\gdn$
in the form
\begin{subequations}
\ber
\sfgupa\T\sfA^{2}\sfgupb=\dab
,
&&\hspace{2.1em}
\sfgupa\T\sfA\sfD\sfA\hspace{0.1em}\sfgupb
=\lambda_{\alpha}\dab
,\\
\sfgdna\T\sfA^{-2}\sfgdnb=\dab
,
&&
\sfgdna\T\sfA^{-1}\sfD\sfA^{-1}\hspace{0.1em}\sfgdnb
=\lambda_{\alpha}\dab
,
\eer
\end{subequations}
and also
\be
\sfgupa\T\sfgdnb=\dab,\qquad\sfgdna\T\sfgupb=\dab
.
\label{updn}
\ee
In this matrix notation we repeat eq.~(\ref{siupa}) as
\be
\label{analcontSG3}
\hsiupa=\sumbpot \left(\sfg_\alpha\T\sfA^{-1}\sfg_\beta\right)
\hsp \hat{s}_\beta
.
\ee
We note for future reference that, although the transformation matrix
$\sfg_\alpha\T\sfA^{-1}\sfg_\beta$ may be banded, it is not in general
possible to truncate it to circumvent having to calculate the full set
of $\hat{s}_\beta$ even if we are only interested in a truncated set
of coefficients $\hsiupa$.

Finally, eq.~(\ref{gortho}) and
eqs~(\ref{bandlgup})--(\ref{analcontSG3}) can now be combined to prove
the equivalent results:
\be
\sumapot s_\alpha g_\alpha=\sumapot\siupa\gdna=\sumapot\sidna\gupa
,
\label{equiv3}
\ee
which we use extensively in subsequent sections.

\section{P~o~t~e~n~t~i~a~l{\hsps}f~i~e~l~d{\hsps}e~s~t~i~m~a~t~i~o~n}
\label{pfe}

We return to solving the geodetic problem stated in
Section~\ref{ps}. We are given noisy data, $d$, taken by a satellite
at an altitude, $a$, over an incomplete sampling domain, $R$, and
attempting to estimate the potential field, $s$, that gives rise to
these observations, at its source level on the unit sphere,
$\Omega$. Although the source field has an infinite bandwidth, we will
practically only be able to make bandlimited estimates of it, which we
denote by $\hat{s}$. The spectral limitation to the bandwidth $L$ as well
as the spatial restriction of the observation domain to the region $R$
motivates our seeking an estimate in terms of a set of basis functions
that are spatiospectrally concentrated, rather than using the
non-localized spherical harmonics $\Ylm$ of more conventional
approaches. This new function set is the Slepian basis, $g_\alpha$,
constructed in Sections~\ref{ssp}--\ref{eotb} in a variety of
geometries, but most notably for the axisymmetric case of a
latitudinal belt around the equator, and its complement the double
polar cap, representative of the polar gap in geodesy.

That the geodetic estimation problem is essentially a problem of
spatiospectral localization can be understood by considering a naive
-- and in practice unsuitable -- estimation scheme. Suppose we
construct estimate in the form of eq.~(\ref{shatdef1}), 
\be
\hat{s}=\sumshortL\hat{s}_{lm}\Ylm
,\label{shatdefn1}
\ee
by minimizing its aggregate squared misfit with the data over the
sphere, given by eq.~(\ref{datdefs}). This amounts to solving the
variational problem
\be
\Phi=\intr(\siuph-d)^2\domg=\mbox{minimum}
,
\label{variational0}
\ee
where the integration domain is the region $R$ in which observations
are available. Substituting eqs~(\ref{datdefs}) and~(\ref{shatdef2})
into eq.~(\ref{variational0}) and requiring the partial derivatives
$\pl\Phi/\pl\hat{s}_{lm}$ to vanish yields the condition
\be
\intr \siuph\Ylm\domg=\intr d\hsp \Ylm\domg
\label{mincon}
,
\ee
while the result $\pl^2\Phi/\pl\hat{s}_{lm}^2>0$ as long as $R\ne 0$
guarantees the convexity of the penalty function $\Phi$. Inserting the
representation~(\ref{shatdef2})--(\ref{sigdef4}) into
eq.~(\ref{mincon}) and using the definition of the localization
kernel~(\ref{Dlmlmpdef}) and its inverse~(\ref{inverse}), the estimate
of the field coefficients at source level is given by:
\be
\hat{s}_{lm}=(1+a)^{l+1}\sumshortLp \Dlmlmp^{-1}\intr d\,\Ylmp\domg
\label{hats}
.
\ee
Thus, the estimate depends on the inverse of the localization kernel
$\sfD$. It is therefore directly influenced by the size and the shape
of the region of missing data, as well as by the chosen
bandwidth. Since $\sfD$ tends to have a very low condition number
(see, e.g., Figure~\ref{sdvals}), finding a stable inverse $\sfD^{-1}$
is problematic: the geodetic inverse problem is ill-conditioned, as is
widely advertised even without reference to the localization nature of
the problem \cite[]{Xu1992a,Xu1992b}.

In the following sections we will derive alternative solutions whose
quality we will judge using standard statistical measures
\cite[e.g.][]{Cox+74,Bendat+2000}. The first will be the average of
the squared difference between a single estimate and the mean of all
estimates over a set of realizations of the data, the estimation variance:
\be
v=\langle(\shat-\langle\shat\rangle)^2\rangle=
\langle\shat^2\rangle-\langle\shat\rangle^2
.
\label{eqvariance}
\ee
The angular brackets denote averaging over the ensemble of repeated
observations, each observation being influenced by a different
realization of the random noise. Similarly, we compute the difference
between the mean of the estimators and the unknown signal, the
estimation bias:
\be
b=\langle\shat\rangle-s
.\label{eqbias}
\ee
We refer to the difference between an estimate and the unknown signal
as the estimation error:
\be
\epsilon=\shat-s
.\label{eqerror}
\ee
Finally, we compute the sum of the variance and the squared bias term,
known as the mean-square error, or mse:
\be
\langle\epsilon^2\rangle=v+b^2
.\label{eqmse}
\ee
For the moment we regard the unknown source signal $s$ as the unique
``truth'', i.e. we consider $s$ to be non-stochastic, although the
data derived from it are contaminated by stochastic noise, see
eqs~(\ref{datdefs})--(\ref{white}). 

\section{S~p~h~e~r~i~c~a~l{\hsps}h~a~r~m~o~n~i~c{\hsps}s~o~l~u~t~i~o~n}
\label{shs}

We have seen that a naive least-squares solution to the geodetic inverse
problem in the spherical harmonic basis yields a solution~(\ref{hats})
that is dependent on the inverse of the localization matrix and
therefore in general impossible to stably
compute. One of the many approaches to circumvent this difficulty is by
adding a model norm to the penalty function
\cite[e.g.][]{Hoerl+70a,Hoerl+70b,Marquardt70,Jackson79};
eq.~(\ref{variational0}) only minimized the norm of the data
misfit. In this section we discuss the solution
to this so-called damped least-squares approach. 

\subsection{Damped least-squares approach}

To stabilize the solution we amend the variational problem of
eq.~(\ref{variational0}) by including a weighted model norm:
\be
\intr(\siuph-d)^2\domg+\eta\intbr\siuph^2\domg=\mbox{minimum}
,
\label{variational}
\ee
where $\eta\ge 0$ is a damping parameter. Retaining the spherical
harmonic basis, once again we supply the bandlimited estimate 
\be
\hat{s}=\sumshortL\hat{s}_{lm}\Ylm
,
\label{SHbasis}
\ee
and minimize~(\ref{variational}) with respect to the unknown
coefficients $\shat_{lm}$. After minimal algebra, involving
eqs~(\ref{shatdef2})--(\ref{sigdef4}), (\ref{Dlmlmpdef}) and
(\ref{blabla})--(\ref{inverse}), we obtain 
the spectral-domain solution,
\ber
\lefteqn{\hat{s}_{lm}=(1+a)^{l+1}\sumshortLp
\left(\Dlmlmp+\eta\bDlmlmp\right)^{-1}}\nnr\\
&&\hspace{7em}\times\intr d\,\Ylmp\domg
\label{hatss}
,
\eer
which only holds at the degrees $l\le L$, since, when $l> L$, no
estimate is available, $\hat{s}_{lm}=0$. The case where $a=0$ and
$\eta=1$, for which, from eq.~(\ref{blabla}), $\sfD+\sfbD=\sfI$, the
identity matrix, was treated in some detail by \cite{Sneeuw+97}. The
integral over the data in eq.~(\ref{hatss}) is made explicit by
substituting eq.~(\ref{datdefs}) and using eqs~(\ref{shatdef2})
and~(\ref{Dlmlmpdef}) once again:
\be
\intr d\,\Ylm\domg=\sumshortp\Dlmlmp\suplmp+
\intr n\,\Ylm\domg\label{threeparts}
.
\ee
Comparing eq.~(\ref{hatss}) to eq.~(\ref{hats}), we now require the
inverse of the weighted sum of the operator localizing to $R$ and the
complementary operator localizing to the region of missing data
$\barr$. The addition of the small quantity $\eta\sfbD$ to the
original matrix $\sfD$ improves its condition number. We postpone a
discussion on determining the ideal value of the weighting parameter
$\eta$ but it is clear that the estimate of the field coefficients
$\hat{s}_{lm}$ in the form of eq.~(\ref{hatss}) is now computable.

In order to ascertain the statistical
properties~(\ref{eqvariance})--(\ref{eqmse}) of the new
estimate~(\ref{hatss})--(\ref{threeparts}) we first calculate the
average of this estimate over all realizations of the noise. From
eq.~(\ref{zmean}), this ensemble averaging of
eqs~(\ref{hatss})--(\ref{threeparts}) annihilates the random noise
term, and we obtain
\ber
\label{solv}
\lefteqn{\langle\shat_{lm}\rangle=(1+a)^{l+1}\sumshortLp
\left(\Dlmlmp+\eta\bDlmlmp\right)^{-1}}\nnr\\
&&{}\hspace{7em}\times\hsp\hsp\sumshortpp\Dlmplmpp\suplmpp
.
\eer
Again, the coefficients $\shat_{lm}$ are defined only in the degree
range $l\le L$. We note that, were the source signal to be similarly
bandlimited, the coefficients $\shat_{lm}$ obtained by undamped
($\eta=0$) estimation would be equal to the true source coefficients
$s_{lm}$. This follows directly from substituting the leftmost term of
eq.~(\ref{sigdef3}) into eq.~(\ref{solv}) and using
eqs~(\ref{slmdef3}) and~(\ref{inverse}). The addition of the damping
term ($\eta>0$) biases the estimate away from the truth, hence the
name ``biased estimation'' for this procedure
\cite[]{Hoerl+70b}. It is the price we pay to be able to calculate the
estimate at all.

There are other benefits as well. These are most easily seen by
computing a spatial-domain representation of the estimate using the
Slepian basis, as in eq.~(\ref{shatdef1}). Making use of the
equivalence~(\ref{equiv3}), we write for the (bandlimited) estimate
\be
\hat{s}=\sumapot\hsiupa \gdna 
,\label{altdef}
\ee
noting that the upward continued coefficients $\hsiupa$ in the Slepian
basis are calculated according to eq.~(\ref{siupa}), and the downward
continued Slepian basis functions according to eq.~(\ref{bandlgdn}). A
Slepian basis expansion of the (broadband) observations, combining
eqs~(\ref{datdefs}) and~(\ref{shatdef2}), is given by
\be
d=\sumapot\siupa
g_\alpha+\sumshortR\suplm\Ylm+n
.\label{dexpa}
\ee
This equation allows us to find an alternative expression for the data
integral~(\ref{threeparts}), for which we also use
eqs~(\ref{inversion}), (\ref{Dlmlmpdef}) and the double orthogonality
of the Slepian functions~(\ref{orthogspace}), namely
\ber
\intr d\,\Ylm\domg&=&\sumapot\glma\left(\lambda_\alpha\siupa
+\intr n\hsp g_\alpha \domg\right)\nnr\\
&&{}+\sumshortRp\Dlmlmp\suplmp
.
\label{threeparts2}
\eer
Inserting eqs~(\ref{siupa}), (\ref{hatss}) and~(\ref{threeparts2})
into eq.~(\ref{altdef}) and using the expressions~(\ref{hlm})
and~(\ref{invfut}) yields the spatial-domain estimate of the field as 
\begin{subequations}
\label{sqav}
\ber
\hat{s}(\rhat)&=&\sumapot\lstar \gdna(\rhat)
\\
&&{}\hspace{-1em}\times\left(\lambda_\alpha\siupa+\intr n\hsp
g_\alpha\domg+\sumshortR \halm\suplm\nnr
\right),\\
\lstar&=&\lstart
.
\eer
\end{subequations}
We have introduced the symbol $\lstar$ for notational convenience. In
the absence of damping, $\lambda^*_\alpha(0)=\lambda_\alpha^{-1}$,
i.e. the inverse of the concentration eigenvalue. Here, too, the
necessity of damping is readily apparent: as the eigenvalues of the
concentration operator, $\lambda_\alpha$, become vanishingly small,
their inverse grows explosively, inflating the noise term and the term
containing the signal at the unmodeled degrees $l>L$, and rendering the
stable computation of the estimate~(\ref{sqav}) impossible. Adding the
damping factor is a useful way to prevent this. 

The zero mean of the  stochastic noise, eq.~(\ref{zmean}), guarantees
that the ensemble average of the spatial estimate over all
realizations of the noise is given by 
\ber
\label{hatssp}
\langle\hat{s}(\rhat)\rangle&=&\sumapot\lstar\gdna(\rhat)
\\
&&{}\hspace{2em}\times\left(\lambda_\alpha\siupa
+\sumshortR\halm\suplm\right)\nnr
.
\eer
Eq.~(\ref{hatssp}) can be combined with eqs~(\ref{equiv3})
and~(\ref{sigdefx}) to show that for bandlimited source fields and in
the absence of damping, the mean of the spatial estimate
$\langle\hat{s}(\rhat)\rangle$ is identical to the source field
$s(\rhat)$, even if the estimate $\hat{s}(\rhat)$ is impossible to
compute stably without the damping term.

The introduction of the damping term stabilizes the solution at the
cost of added bias. Following eq.~(\ref{eqbias}) the latter is
calculated by subtracting the full representation of the
signal~(\ref{sigdefx}) from eq.~(\ref{hatssp}), making use of
eqs~(\ref{invfut}) and~(\ref{equiv3}). The spatial estimation bias is
then given by
\ber
\lefteqn{b(\rhat)=-\eta\sumapot(1-\lambda_\alpha)\lstar s_\alpha
g_\alpha(\rhat)}\nnr
\\
&&{}\hspace{3em}+
\sumapot\lstar\gdna(\rhat)\sumshortR\halm\suplm\nnr\\
&&{}\hspace{3em}-\sumshortR s_{lm}\Ylm(\rhat)\label{biaseq}
.
\eer
We have brought forward the damping parameter $\eta$ by using
the identity $\lambda^*\lambda-1=-\eta(1-\lambda)\lambda^*$. In the
absence of damping ($\eta=0$), the first term in this equation
vanishes, leaving us with the unavoidable broadband leakage (the
second term) and bias due to making bandlimited estimates of broadband
fields (the third term). 

An expression for the estimation variance from~(\ref{eqvariance}) is
obtained by squaring eq.~(\ref{sqav}) and averaging the result, using
the properties of the noise~(\ref{zmean})--(\ref{white}) and
eq.~(\ref{orthogspace}), and subtracting from the result the square of
eq.~(\ref{hatssp}). The spatial estimation variance is
\be
v(\rhat)=N\sumapot\lambda_\alpha[\lstar]^2\gdna^2(\rhat)
\label{vargro}
.
\ee
We note that eq.~(\ref{vargro}) is the only one thus far to assume
that the power spectrum of the noise is white, of magnitude $N$. And
one more time the necessity of the damping is apparent: in its
absence, the estimation variance strongly amplifies the measurement
noise. At the price of introducing additional bias, damping prevents
this.

\subsection{A bandlimited white stochastic source}

In the previous section we have derived expressions for the average
estimate of the spherical harmonic field coefficients,
$\langle\shat_{lm}\rangle$, in eq.~(\ref{solv}), and for the average
of spatial expansions of the estimated field, $\langle
s(\rhat)\rangle$, in eq.~(\ref{hatssp}). The averaging was over the
different realizations of the stochastic noise process. Both
expressions are valid in the most general sense; the only condition
being that the average over all realizations of the noise, $\langle
n(\rhat)\rangle$, is zero. No further assumptions are necessary. We
have drawn attention to the fact that without the damping term, both
estimates are nearly impossible to calculate. However, in that case,
they are unbiased when the source signal itself is strictly
bandlimited to within a bandwidth $L$ identical to that of the
estimate.

We can make this explicit by postulating that the geophysical signal
expressed as eq.~(\ref{sigdef2}) or eq.~(\ref{sigdefx}) has spherical
harmonic expansion coefficients that vanish outside of this bandwidth:
\be
s_{lm}= 0\quad\mbox{for}\quad L< l\le \infty
.
\label{disreg}
\ee
We will work with this contrived geophysical signal for the simple
reason that no amount of sophistication can cure the fact that forming
harmonically truncated estimates leads to multiple bias terms, as can
be seen from eq.~(\ref{biaseq}). Under the condition~(\ref{disreg}),
eq.~(\ref{solv}) becomes
\ber
\lefteqn{\langle\shat_{lm}\rangle=(1+a)^{l+1}\sumshortLp
\left(\Dlmlmp+\eta\bDlmlmp\right)^{-1}}\nnr\\
&&{}\hspace{1.5em}\times\hsp\hsp\sumshortLpp\Dlmplmpp\suplmpp
,\label{newsolve}
\eer
from which, using eqs~(\ref{slmdef3}) and~(\ref{inverse}), we derive
immediately that the undamped estimate of the
coefficients, given by eq.~(\ref{hats}), is unbiased: 
\be
\langle\shat_{lm}\rangle=s_{lm}\quad\mbox{if}\quad \eta=0
.
\label{newsolveunbias}
\ee
Similarly, using eq.~(\ref{equiv3}), eq.~(\ref{hatssp}) can be
transformed under the same condition~(\ref{disreg})
into 
\be
\langle\hat{s}(\rhat)\rangle=\sumapot\lambda_\alpha\lstar\hsp s_\alpha
g_\alpha(\rhat) 
\label{newsolveb}
,
\ee
from which, with $\lambda^*_\alpha(0)=\lambda_\alpha^{-1}$, the
undamped spatial estimate of the field, given by eq.~(\ref{shatdefn1}),
is unbiased:
\be
\langle\hat{s}(\rhat)\rangle=s(\rhat)\quad\mbox{if}\quad \eta=0
.
\label{abcd}
\ee
Indeed, under the condition~(\ref{disreg}), the only term left in the
bias equation~(\ref{biaseq}) is directly, though not linearly,
dependent on the damping term $\eta$: it is
\ber
b(\rhat)=-\eta\sumapot(1-\lambda_\alpha)\lstar\hsp s_\alpha\garh
\label{newsolvebb}.
\eer

Although we can calculate the mean-square estimation
error~(\ref{eqmse}) exactly from eqs~(\ref{vargro}) and\
(\ref{newsolvebb}), we will gain additional insight when we cease to
consider the unknown signal as a non-stochastic signal.  The source
signal $s(\rhat)$, until now, has been considered to be ``given'': we
have simply assumed it is of the form~(\ref{sigdef2}) and attempted to
estimate its true unknown coefficients $s_{lm}$ from incomplete and
noisy observations. All averaging in the construction of the bias and
variance terms was carried out over the different realizations of the
noise $n(\rhat)$, which we took to be a white stochastic process. By
now considering the geophysical signal, as well, to be a stochastic
process, we shall calculate the mse after an additional round of
averaging, this time over the various realizations of $s(\rhat)$,
should they be available. Instead of eq.~(\ref{eqmse}) we thus 
write
\be
\quad\langle\epsilon^2\rangle=v+\langle b^2\rangle
,
\label{nepsi}
\ee
where the angular brackets now denote an average over the ensemble of
signals. Strictly speaking we should write
$\langle\langle\epsilon^2\rangle\rangle$ but we eschew the double
brackets in the interest of notational simplicity.

We notice from eq.~(\ref{newsolvebb}) that to compute
$\langle b^2\rangle$ we shall require the covariance $\covsa$ of the
expansion coefficients of the field in the Slepian basis. To
facilitate the treatment and for easy comparison with the assumed
white power spectrum of the noise process, we shall consider a
bandlimited source signal that is ``whitish'', i.e. white within the
band $l\le L$, such that its covariances in the spherical harmonic and
Slepian bases, respectively, are given by
\begin{subequations}
\label{covSG}
\ber
\covsl&=&S\hsp\dllp\dmmp
,\\
\covsa&=&S\hsp\dab
,
\eer
\end{subequations}
while noting that, as far as the spatial covariance of this signal
concerned, 
\be
\label{noteq}
\covsr=S\,D(\rhat,\rhat')\ne S\,\delta(\rhat,\rhat')
,
\ee
as can be deduced by combining eq.~(\ref{sigdef2}) with
eq.~(\ref{covSG}) and using eqs~(\ref{additionSH})--(\ref{dirac})
and~(\ref{banddelta}). A last assumption introduced here is that the
noise is wholly uncorrelated with the signal:
\be
\langle s(\rhat)n(\rhat')\rangle=0.\label{uncor}
\ee
The average of the squared bias term~(\ref{newsolvebb}) under these
idealized assumptions is
\be
\langle b^2(\rhat)\rangle=\eta^2 S\sumapot
(1-\lambda_\alpha)^2[\lstar]^2 g_\alpha^2(\rhat)
,\label{biasgro}
\ee
and the mean-square estimation error, following eq.~(\ref{nepsi}), is
formed by combining this result with the expression for the variance
in eq.~(\ref{vargro}). The latter expression is unchanged even if the
source signal is stochastic, as long as the noise is uncorrelated with
the signal, eq.~(\ref{uncor}). Thus, the mean-square error of the
bandlimited estimation of a bandlimited white source field from
incomplete observations at an altitude in the presence of white noise
is given by
\ber
\label{msefinal}
\langle\epsilon^2(\rhat)\rangle&=&N\sumapot\lambda_\alpha[\lstar]^2
\gdna^2(\rhat)\\
&&{}+\eta^2 S\sumapot
(1-\lambda_\alpha)^2[\lstar]^2 g_\alpha^2(\rhat)\nnr
.
\eer
All $\Lpot$ basis functions are required to form the mse. The first
term in the expression for the mse is the variance: it is the only
term that depends on the noise. We have seen that without damping
($\eta=0$) this term becomes unmanageably large: the addition of
damping counteracts this. In addition, the estimation variance also
varies with the observation height $a$ above the unit sphere: as $a$
grows, so do the downward continued Slepian basis functions
$\gdn(\rhat)$, and with them, the noise. The second term in the mse is
due to bias. This is the only term that depends on the characteristics
of the signal. It is independent of the satellite altitude at which
the measurements are taken.

\begin{figure*}\center
\rotatebox{0}{
\iftwocol{
\includegraphics[width=0.7\textwidth]{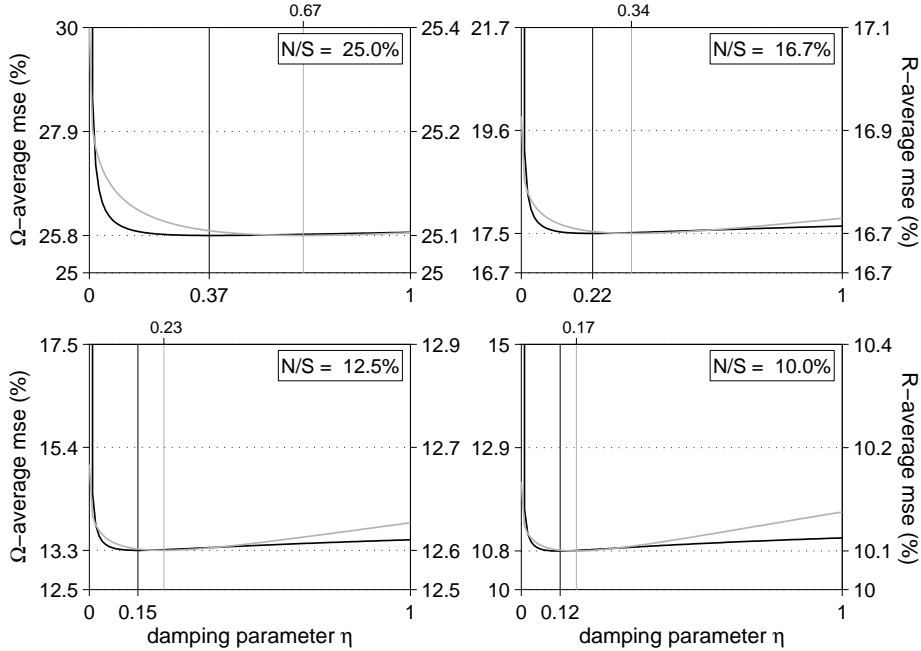}}
{\includegraphics[width=0.85\textwidth]{sdeta.eps}}
}  
\caption{Spatially averaged mean-square error (as a percentage of the
average mean-square signal strength) of the damped-least-squares
spherical harmonic solution to the geodetic estimation problem for a
bandlimited white signal and white noise. The $\Omega$-average mse
(black curves and left ordinate) is the average over the entire
sphere; the $R$-average mse (grey curves and right ordinate) is the
average over the region of observation, the equatorial belt
complementary to the polar gap of radius $\Theta=10^{\circ}$. The
bandwidth of the signal and its estimate is $L=45$. The measurement
altitude is $a=0$. The signal-to-noise levels shown are $S/N= 4, 6, 8$
and $10$. We plot the normalized average mse values as a function of
the damping parameter $\eta$, and indicate by vertical lines the
values of $\eta$ that minimize them.  The range of the $R$-average is
much reduced compared to the $\Omega$-average values. Both ordinates
are truncated below at $N/S$, the mse value when the observation
region is the entire sphere, $R=\Omega$.}
\label{sdeta}  
\end{figure*} 

\begin{figure*}\center
\rotatebox{0}{\iftwocol{
\includegraphics[width=0.7\textwidth]{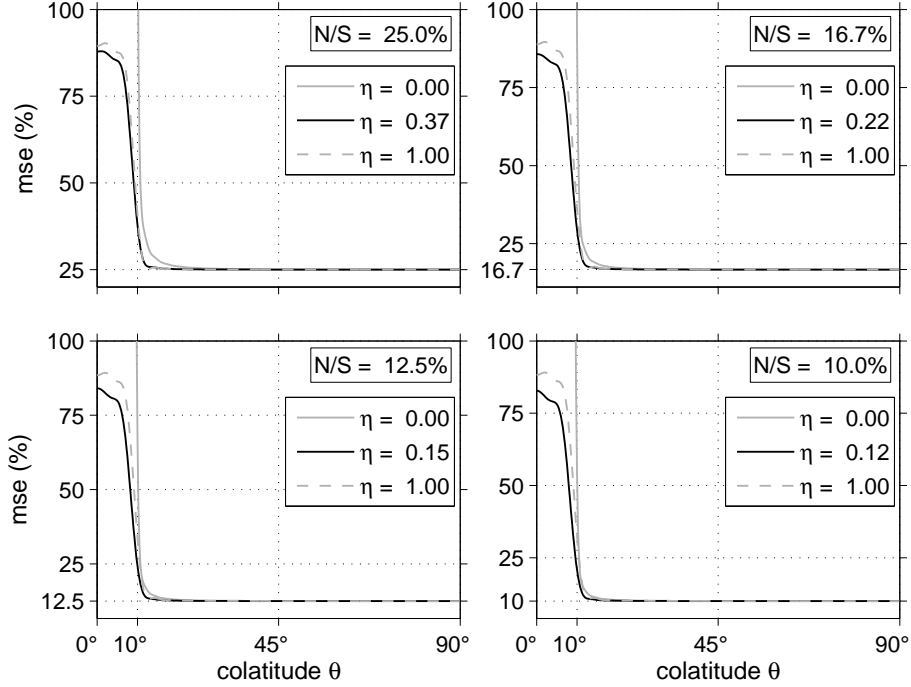}
}{
\includegraphics[width=0.85\textwidth]{sderrs.eps}}
} 
\caption{Mean-square error (as a percentage of the mean-square
signal strength) of the damped-least-squares spherical harmonic
solution to the geodetic estimation problem for a bandlimited white
signal and white noise. As in Figure~\ref{sdeta}, the signal-to-noise
levels shown are $S/N= 4, 6, 8$ and $10$, the measurement altitude
$a=0$, the polar gap consists of caps with radius $\Theta=10^{\circ}$,
and the bandwidth of the signal and its estimate is $L=45$. As a
function of colatitude, for an arbitrary longitude, we plot the mse of
the undamped solution ($\eta=0$), of a heavily damped solution
($\eta=1$), and of the solution at the damping level which minimizes
the normalized average mse over the unit sphere, i.e. for the values
$\eta_\Omega=0.37, 0.22, 0.15$, and $0.12$, that were marked by black
vertical lines in Figure~\ref{sdeta}. The mse is symmetric about the
equator. The ordinate is truncated at 100\%; the mse of the undamped
solution in the region of the polar gap exceeds this value by several
orders of magnitude. }
\label{sderrs} 
\end{figure*} 

\subsection{Optimal damping level}

To illustrate the behavior of the mse in eq.~(\ref{msefinal}) we will
focus on the case where the measurement altitude is $a=0$, hence
$g_\alpha=\gdna=\gupa$. This simplifies the expressions to:
\begin{subequations}
\label{msefin2}
\ber
\langle\epsilon^2(\rhat)\rangle&=&\sumapot
\ffae g_\alpha^2(\rhat),\\
\ffae&=&[\lstar]^2
\left[N\lambda_\alpha+\eta^2S(1-\lambda_\alpha)^2\right].
\eer
\end{subequations}
The function $\ffae$ combines the effects of data noise, damping,
signal strength, and measurement geometry. We will compare the
mean-square error with the mean-square signal strength over all
realizations, which is given by
\be
\langle s^2(\rhat)\rangle=S\hsp\frac{\Lpot}{4\pi}
.
\label{slpot}
\ee
The result~(\ref{slpot}) is obtained by combining eq.~(\ref{noteq})
with the definition~(\ref{banddelta}) at $\rhat\cdot\rhat=1$. We
calculate the following two quantities. First, a normalized spatial
average of the mse given by the ratio of the mean square
error~(\ref{msefin2}) to the mean square signal
strength~(\ref{slpot}), both averaged over the entire sphere
$\Omega$. Using the orthogonality conditions~(\ref{orthogspace}) this
``$\Omega$-average mse'' is given by
\be
\fracd{\into\langle\epsilon^2(\rhat)\rangle\domg}
{\into\langle s^2(\rhat)\rangle\domg}=
\frac{1}{S\Lpot}\sumapot {\ffae}
.
\label{mseO}
\ee
Second, a scaled ``$R$-average mse'' is given by the ratio of the same
quantities, averaged over the covered region $R$. Using
eq.~(\ref{orthogspace}) and the definition~(\ref{tracedef}) of the
Shannon number $K$, it is
\be
\fracd{\intr\langle\epsilon^2(\rhat)\rangle\domg}
{\intr\langle
 s^2(\rhat)\rangle\domg}=
\frac{1}{KS}\sumapot {\lambda_\alpha\ffae}
.\label{mseR}
\ee
Both quantities are shown in Figure~\ref{sdeta}, for a double-cap
polar gap of $\Theta=10^{\circ}$ and a bandwidth $L=45$. They are
plotted in different panels for different signal-to-noise ratios
$S/N=4,6,8$ and $10$ as functions of the damping parameter $\eta=0\rar
1$. We show eq.~(\ref{mseO}) in black, with the scale on the left of
the panels, and eq.~(\ref{mseR}) in grey, with the scale on the right
hand side. The range of $\Omega$-average mse values shown is much
larger (5\% in all four panels) than the equivalent range in
$R$-average mse values (0.4\% in all panels): the effects of damping
on the overall mse over the entire globe are much more pronounced than
its effects on the mse averaged over the region in which data were
collected. The ordinate is truncated to aid the visualization. The
maximum $R$-average mse is $(N/S)/\cos\Theta$ which is attained when
$\eta=0$. This can be verified by noting that $\ffan=N\lambda_a^{-1}$,
using the definition of the Shannon number~(\ref{tracedef}), and
noting that the area of the covered region is equal to
$A=4\pi\cos\Theta$. Thus, at a given signal-to-noise ratio only the
size of the polar gap controls the upper bound on the
$R$-average mse. A lower bound for all damping levels is found at full
coverage, $R=\Omega$. It thus applies to both measures of the average
mse. Indeed without a polar gap, $\Theta=0$, $K=\Lpot$,
$\lstar=\lambda_{\alpha}=1$, $\ffae=N$, and the scaled average mse
curves never drop below $N/S$, which we use as a lower cutoff for the
vertical axes.

A statistically desirable estimator \cite[e.g.][]{Cox+74,Bendat+2000}
is one that is unbiased and efficient, i.e. it minimizes the mean
square estimation error. We have seen that sacrificing the
unbiasedness by introducing damping removes the obstacles in computing
the estimate in the first place, and reduces the estimation variance.
We can calculate the damping level that is overall optimal by
minimizing the mse~(\ref{msefin2}) with respect to the damping
parameter $\eta$. However, minimization of the $R$-average and
$\Omega$-average mse will yield slightly different optima. Minimizing,
eq.~(\ref{mseO}), the normalized mse over the entire sphere we obtain
an optimal damping coefficient $\eta_\Omega$ given by
\be
\eta_\Omega=\frac{N}{S}\frac{\sumapot[\lstarO]^3\lambda_\alpha
(1-\lambda_\alpha)}  
{\hspace{0.6em}\sumapot[\lstarO]^3\lambda_\alpha(1-\lambda_\alpha)^2}
\label{eO}
.
\ee
Likewise, minimization of eq.~(\ref{mseR}), the normalized mse over
the region of coverage, yields an optimal damping coefficient $\eta_R$
given by 
\be
\eta_R=\frac{N}{S}\frac{\sumapot[\lstarR]^3\lambda_\alpha^2
(1-\lambda_\alpha)}
{\hspace{0.6em}\sumapot[\lstarR]^3\lambda_\alpha^2(1-\lambda_\alpha)^2}
\label{eR}
.
\ee
Although the unknown optimal damping levels $\eta_\Omega$ and $\eta_R$
appear on both sides of eqs~(\ref{eO}) and~(\ref{eR}), their values
can be easily computed by iteration. They depend on the measurement
geometry, the damping, and the signal-to-noise ratio. In
Figure~\ref{sdeta}, $\eta_\Omega$ and $\eta_R$ are shown as black and
grey vertical lines, respectively. At high signal-to-noise ratios both
can be approximated as $\eta_R=\eta_\Omega\approx N/S\ll 1$.

When the coverage region is axisymmetric the mse~(\ref{msefin2}) is
independent of the longitude, as can be deduced from
eq.~(\ref{polarg2}). Thus, in Figure~\ref{sderrs} we plot
$\langle\epsilon^2(\theta)\rangle/\langle s^2(\theta)\rangle$, in
percent, for different signal-to-noise ratios, as a function of
colatitude and for various damping levels: i.e. in the undamped
($\eta=0$), fully damped ($\eta=1$) and optimally damped case
($\eta=\eta_\Omega$). The vertical axes are truncated at 100\% as the
undamped values exceed this value by many orders of magnitude.

\section{S~l~e~p~i~a~n{\hsps}b~a~s~i~s{\hsps}s~o~l~u~t~i~o~n}
\label{sbs}

In the previous section we expanded the estimate of the signal into a
bandlimited spherical harmonic basis and performed a damped
least-squares inversion for the unknown coefficients. This estimation
procedure resulted in a biased estimate, but the damping prevented the
detrimental amplification of the measurement noise. We derived
expressions for the optimal level of damping required for ``whitish''
signals measured at zero altitude. Adding a small amount of bias made
the estimate computable and reduced its variance. We used the
(downward continued) Slepian basis to find expressions for the
resultant damped spherical harmonic estimate in the spatial domain and
to find its bias, variance, and mse. Using the Slepian basis greatly
simplified the expressions because of the fact that, as opposed to the
spherical harmonics, the Slepian functions are orthogonal over both
the entire sphere and the closed subdomains over which, by design,
their energy is optimally concentrated.

Alternatively, we might have sought an estimate that is expressed in
the spherical Slepian basis at the start. As we have seen, the first
$K$ Slepian eigenfunctions, where $K$ is the Shannon
number~(\ref{tracedef}), provide an excellent coverage of the region
of observation. This implies that their associated eigenvalues
$\lambda$ are close to unity, avoiding any problems with their
inversion. In this section we will explore the effect on the geodetic
solution of using a truncated Slepian basis, consisting of the
$J$ basis functions that are best concentrated over the region of
satellite observation. Even if $J=K$ appears to be a natural choice,
we will determine the truncation level $J$ by optimization of the
mean-square estimation error, as we did to find the optimal damping
parameter in the damped least-squares spherical harmonic approach.

\subsection{Truncated Slepian function approach} 

The original undamped problem posed in eq.~(\ref{variational0}),
\be
\intr(\siuph-d)^2\domg=\mbox{minimum}
,\label{variational2}
\ee
is now solved by expanding the estimate in the downward continued
truncated Slepian basis
\be
\hat{s}=\sumak\hsiupa \gdna
,
\label{SGest1}
\ee
and minimizing eq.~(\ref{variational2}) with respect to the estimation
coefficient $\hsiupa$. The second derivative of
eq.~(\ref{variational2}) is always positive. After minimal
algebra, using eq.~(\ref{shatdef2}) and the double
orthogonality~(\ref{orthogspace}), the expansion coefficients in
eq.~(\ref{SGest1}) are obtained from
\be
\hsiupa=\lambda_\alpha^{-1}\intr d\hsp g_\alpha\domg
.\label{SGsolspec}
\ee
This result can alternatively be derived by substituting
eq.~(\ref{hatss}) of the damped spherical harmonic approach into
eq.~(\ref{siupa}), setting $\eta=0$, and using eqs~(\ref{inversion})
and~(\ref{invfut}). We purposely chose an estimate the
form~(\ref{SGest1}) to find the truncated expansion coefficients
$\hsiupa$ in their upward continued form and multiplying the downward
continued Slepian functions $\gdna$, rather than simply expressing
eq.~(\ref{hats}) in the Slepian basis $g_\alpha$. In the latter case,
as can be readily verified by combining eq.~(\ref{hats}) with
eqs~(\ref{rota}), (\ref{inversion}) and~(\ref{invfut}), every one of
the expansion coefficients $\hat{s}_\alpha$ would depend on a linear
combination of all $\Lpot$ terms $\lambda_\beta^{-1}\int_Rd\hsp
g_\beta\domg$ through a matrix term $\sfg_\alpha\T\sfA\,\sfg_\beta$
whose kind we have encountered in eq.~(\ref{analcontSG3}). This would
therefore invalidate the method of truncation as a means to avoid the
difficult-to-compute and unnecessarily influential large inverse
eigenvalues. By choosing the representation~(\ref{SGest1}) instead, we
take advantage of eq.~(\ref{equiv3}) to juxtapose upward and downward
continuation, $(1+a)^{l+1}(1+a)^{-l-1}$, thereby canceling their
effect altogether: eq.~(\ref{SGsolspec}) shows that every coefficient
$\hsiupa$ only depends on the inverse eigenvalue at the same rank
$\alpha$. The effect of the measurement at altitude has not
disappeared: it is now contained in eq.~(\ref{SGest1}) in the basis
$\gdna$, of which only the first $J$ functions are required. These are
calculated via eq.~(\ref{bandlgdn}) and ultimately, by the stable
Gr\"unbaum algorithm central to our analysis.

The data integral~(\ref{SGsolspec}) can be calculated by 
substituting into it eqs~(\ref{dexpa}), (\ref{orthogspace}),
(\ref{inversion}), (\ref{Dlmlmpdef}) and (\ref{hlm}), to yield
\be
\intr d\hsp g_\alpha\domg=\lambda_\alpha\siupa+\intr n\hsp
g_\alpha\domg+\sumshortR 
\halm\suplm
.
\label{SGdint}
\ee
Averaging the expressions~(\ref{SGsolspec})--(\ref{SGdint}) over many
estimates annihilates the influence of the random noise by virtue of 
eq.~(\ref{zmean}), and gives
\be
\langle\hsiupa\rangle=\siupa+\lambda_\alpha^{-1}\sumshortR\halm\suplm 
.
\label{bloblo}
\ee
Combining eqs~(\ref{SGest1})--(\ref{SGdint}) yields the estimate in
the space domain,
\ber
\hat{s}(\rhat)&=&\sumak\lambda_\alpha^{-1}\gdna(\rhat)\label{SGsqav}
\\
&&{}\hspace{0em}\times\left(\lambda_\alpha\siupa+\intr n\hsp
g_\alpha\domg+\sumshortR\halm\suplm\nnr
\right)
,
\eer
which, reassuringly, amounts to the truncated but undamped ($\eta=0$)
version of eq.~(\ref{sqav}). As before we can eliminate the noise term
by averaging over many realizations, to obtain
\ber
\langle\hat{s}(\rhat)\rangle&=&
\sumak\gdna(\rhat)\nnr   
\\
&&{}\hspace{0em}\times\left(\siupa+\lambda_\alpha^{-1}\sumshortR
\halm\suplm\right)\label{SGsqavn}
.
\eer
The estimation bias, following eq.~(\ref{eqbias}), is obtained by
subtracting from eq.~(\ref{SGsqavn}) the representation of the
signal~(\ref{sigdefx}) and using the equivalence~(\ref{equiv3}):
\ber
\lefteqn{b(\rhat)=-\sumakR  s_\alpha
g_\alpha(\rhat)-\sumshortR s_{lm}\Ylm(\rhat)}\nnr 
\label{SGbiaseq} 
\\
&&{}\hspace{3em}+
\sumak\lambda_\alpha^{-1}\gdna(\rhat)\sumshortR\halm \suplm
.
\eer
Without truncation of the Slepian basis function set, i.e. when
$J=\Lpot$, the first term in this equation vanishes. The remaining
contributions arise due to forming bandlimited estimates of broadband
signals, leading to unavoidable leakage and broadband bias. Comparing
eqs~(\ref{biaseq}) and~(\ref{SGbiaseq}) we discover the parallel roles
of damping and truncation. The introduction of the damping parameter
$\eta$ adds an extra bias term to the expression~(\ref{biaseq}), and
reduces the size of the leakage term by which the coefficients
$\halm$, $lm>L$ of eq.~(\ref{normratio2}) make the influence of the
signal outside the bandwidth felt, but it is powerless against the
bias due to the bandlimited approximation of the broadband signal,
which is simply that portion of the signal that is outside the
bandwidth $L$. Similarly, increasing the Slepian truncation level by
the reduction of $J$ from $\Lpot$ in eq.~(\ref{SGbiaseq}) introduces a
new term in the expression for the estimation bias, and reduces the
effect of the leakage term containing the coefficients $h$, but it is
again no match for the remaining broadband bias from the
bandlimitation of the estimate.

An expression for the estimation variance,
eq.~(\ref{eqvariance}), is obtained by squaring and averaging
eq.~(\ref{SGsqav}), using the noise
properties~(\ref{zmean})--(\ref{white}) and the orthogonality of the
Slepian basis functions~(\ref{orthogspace}), and subtracting the
square of~(\ref{SGsqavn}). The resulting variance is
\be
v(\rhat)=N\sumak\lambda_\alpha^{-1}\gdna^2(\rhat)
\label{SGvargro}
.
\ee
This expression is again the first in this section in which we have
used the white noise assumption, and once again it will be valid even
if the source signal is considered stochastic as long as
eq.~(\ref{uncor}) holds. Comparison of the variance expression in this
truncated Slepian basis approach with eq.~(\ref{vargro}) obtained via
the damped spherical harmonics method validates our approach. Without
damping, when $\eta=0$ in eq.~(\ref{vargro}), or without truncation,
$J=\Lpot$ in eq.~(\ref{SGvargro}), both expressions are
identical. Much like the damping term, the truncation of the basis set
to its first $J\le\Lpot$ elements reduces the estimation variance by
checking the growth of the terms $\lambda_\alpha^{-1}$. The more
severe the truncation, the lower $J$, and the lower the variance
becomes.

\subsection{A bandlimited white stochastic source}

Once again, we now focus on geophysical signals that are white within
a bandwidth $L$ as expressed by eqs~(\ref{disreg})
and~(\ref{covSG}). This assumption transforms eq.~(\ref{bloblo}) into
\be
\langle\hsiupa\rangle=\siupa
\label{SGunbiased}
,
\ee
illustrating the fact that an estimate of the form~(\ref{SGsolspec})
is spectrally unbiased. Just as our analysis of the damped spherical
harmonic method showed that for bandlimited source fields, the
undamped, i.e. $\eta=0$, estimate of eq.~(\ref{hats}) is incomputable
due to the ill-conditioning of $\sfD^{-1}$, but unbiased, as shown by
eq.~(\ref{newsolveunbias}), we have now shown that the untruncated,
i.e. $\alpha=1\rar\Lpot$, Slepian basis estimate of
eq.~(\ref{SGsolspec}) is incomputable due to the growth of the
eigenvalues $\lambda^{-1}$, although it, too, is unbiased, as shown by
eq.~(\ref{SGunbiased}). The damping term makes the estimate computable
but biased, just as the truncation of the eigenvalues prevents the
blow-up of their inverse at the cost of added bias.

In the spatial domain, using eqs~(\ref{SGsqavn}) and~(\ref{equiv3}),
the average over all estimates is then
\be
\langle\hat{s}(\rhat)\rangle=
\sumak g_\alpha(\rhat)s_\alpha
.\label{SGsqavnun}  
\ee
In the absence of truncation, $J=\Lpot$, the spatial estimate of the
form~(\ref{SGest1}) is similarly unbiased:
\be
\langle\hat{s}(\rhat)\rangle=s(\rhat)
,
\ee
which we may again compare to the unbiasedness~(\ref{abcd}) of the
undamped estimate~(\ref{shatdefn1}).  
Explicitly, under the condition (\ref{disreg}), the only contributing
term in eq.~(\ref{SGbiaseq}) is given by
\be
b(\rhat)=-\sumakR  s_\alpha
g_\alpha(\rhat)\label{SGbiasequ} 
.
\ee
This term decreases with increasing $J$, and vanishes altogether when
$J=\Lpot$. It can be compared to eq.~(\ref{newsolvebb}). The average
over all realizations of the signal of the squared bias, for a
``whitish'' signal with covariance~(\ref{covSG}), is given by
\be
\langle b^2(\rhat)\rangle=S\sumakR 
 g_\alpha^2(\rhat)
,\label{SGbiasgro}
\ee
which should be compared with the corresponding eq.~(\ref{biasgro}) in
the damped spherical harmonic case.  From this and
eq.~(\ref{SGvargro}) we can calculate the mean-square estimation error
following eq.~(\ref{nepsi}), which is now
\be
\langle\epsilon^2(\rhat)\rangle=N\sumak\lambda_\alpha^{-1}\gdna^2(\rhat)+
S\sumakR g_\alpha^2(\rhat)
.
\label{SGmsefinal}
\ee
The mse of the untruncated Slepian basis approach and that of the
undamped spherical harmonic estimation method~(\ref{msefinal}) are
identical. This of course is a direct consequence of the fact that
both bases are related to each other by the orthonormal transformation
eqs~(\ref{inversion})--(\ref{gortho}). Of note is the very different
form of the damped and truncated expressions, eqs~(\ref{msefinal})
and~(\ref{SGmsefinal}), for the mse. Whereas eq.~(\ref{msefinal})
consists of a weighted sum of all basis functions $\alpha=1\rar\Lpot$
in a manner that appears to mix the influence of the noise, the
damping, and the signal, the truncated expression~(\ref{SGmsefinal})
has disentangled the effects of the noise and the signal by
distributing the influence of the variance over the basis functions
$\alpha=1\rar J$ that are well concentrated inside the measurement
area, and the effect of the bias over those $\alpha=J+1\rar\Lpot$ that
are concentrated in the region of missing data. To the one piece that
is missing, the decision on where to truncate the data by the choice
of $J$, we now turn.

\begin{figure*}\center
\rotatebox{0}{
\iftwocol
{\includegraphics[width=0.7\textwidth]{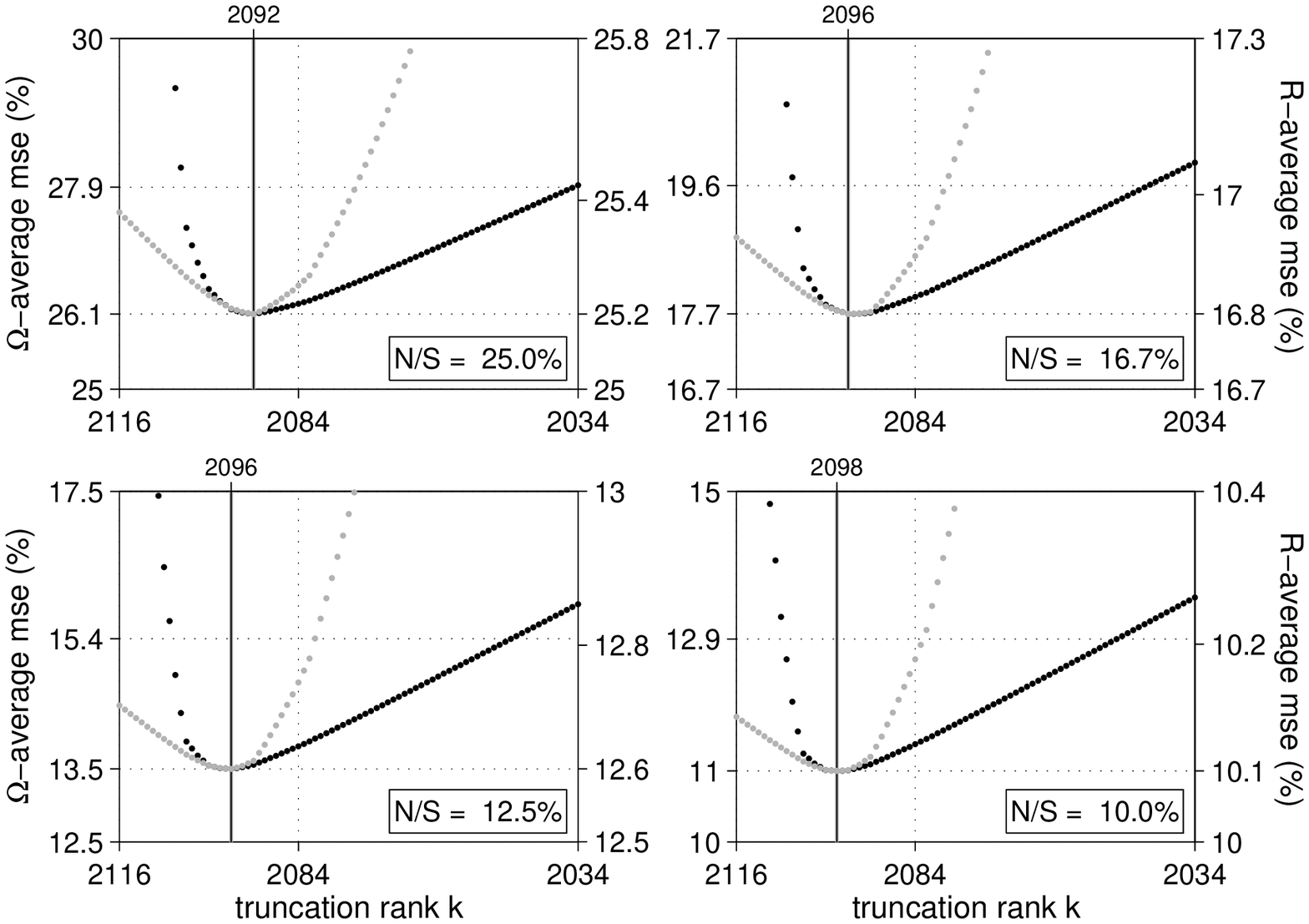}}
{\includegraphics[width=0.85\textwidth]{sdkay.eps}}
} 
\caption{Spatially averaged mean-square error (as a percentage of the
average mean-square signal strength), of the truncated Slepian
function solution to the geodetic estimation problem for a bandlimited
white signal and white noise.  The $\Omega$-average mse (black curves
and left ordinate) is the average over the entire sphere; the
$R$-average mse (grey curves and right ordinate) is the average over
the region of observation, the equatorial belt complementary to the
polar gap of radius $\Theta=10^{\circ}$. The bandwidth of the signal
and its estimate is $L=45$. The measurement altitude is $a=0$. The
signal-to-noise levels shown are $S/N= 4, 6, 8$ and $10$. We plot the
average mse values as a function of the truncation rank $J$, and
indicate by solid vertical lines the values $J_\Omega=J_R$ that
minimize them, and the Shannon number $K=2084$ by the dotted vertical
line. The range of the $R$-average is much reduced compared to the
$\Omega$-average values. Both ordinates are truncated below at $N/S$,
the value when the observation region is the entire sphere,
$R=\Omega$. The abscissa shows truncation levels from $J=K-50$ up to
the untruncated case, $J=\Lpot$. In that case the mse values are
identical to those of the undamped ($\eta=0$) spherical harmonic case
shown in Figure~\ref{sdeta}.}
\label{sdkay}  
\end{figure*}

\subsection{Optimal truncation level}

To illustrate the behavior of the mse in eq.~(\ref{SGmsefinal}) we
again focus on the zero-altitude case, for which
\be
\langle\epsilon^2(\rhat)\rangle=
N\sumak\lambda_\alpha^{-1}g_\alpha^2(\rhat)+ 
S\sumakR g_\alpha^2(\rhat)
.\label{SGmsefinala0}
\ee
In order to find the optimal truncation level, we consider the
full-sphere and coverage-domain average mse (normalized by the
corresponding quadratic signal averages) as in the damped spherical
harmonic approach. Using eqs~(\ref{slpot}) and~(\ref{orthogspace}) we
find from eq.~(\ref{SGmsefinala0}) that the $\Omega$-average mse in
the truncated Slepian case is
\be
\fracd{\into\langle\epsilon^2(\rhat)\rangle\domg}
{\into\langle s^2(\rhat)\rangle\domg}=
1-\fracd{J}{\Lpot} 
+
\fracd{N}{S}\sumak\fracd{\lambda_\alpha^{-1}}{\Lpot}
.\label{SGmseO}
\ee
Using the definition~(\ref{tracedef}) of the Shannon number, we
likewise find the $R$-average mse:
\be
\fracd{\intr\langle\epsilon^2(\rhat)\rangle\domg}
{\intr\langle s^2(\rhat)\rangle\domg}=
\fracd{N}{S}\frac{J}{K}+
\fracd{1}{K}\sumakR {\lambda_\alpha}
.\label{SGmseR}
\ee
Both quantities are plotted in Figure~\ref{sdkay} for different
signal-to-noise ratios $S/N=4,6,8$ and $10$ and with the other
parameters unchanged from those of Figure~\ref{sdeta}: a double-cap
polar gap of $\Theta=10^{\circ}$, a bandwidth $L=45$. In black, with
the scale on the left of the panel, we show eq.~(\ref{SGmseO}) as a
function of the truncation level ranging over $J=K-50\rar
\Lpot$. The abscissa is inverted since $J=\Lpot$ corresponds to a
situation without Slepian truncation, and as $J$ decreases, the degree
of truncation increases. In grey, we plot eq.~(\ref{SGmseR}), with a
much reduced scale on the right hand side. The range of
$\Omega$-average mse values shown is much larger (5\% in all four
panels) than the equivalent range in $R$-average mse values (which
varies from panel to panel but is less than 0.8\%): the effects of
truncation on the overall mse over the entire globe are much more
outspoken than its effects on the mse averaged over the region in
which data were collected. This behavior mimics the one seen in
Figure~\ref{sdeta} for the damped spherical harmonic case. The
ordinate is again truncated for clarity. The value of the untruncated
$R$-average mse, attained when $J=\Lpot$, is $(N/S)/\cos\Theta$. This
follows from the definition of the Shannon number~(\ref{tracedef}) and
the area of the covered region, $A=4\pi\cos\Theta$, and is identical
to the corresponding value in the undamped spherical harmonic case. A
lower bound is found at full coverage, when $R=\Omega$,
$\Theta=0^{\circ}$, $K=\Lpot$ and $\lambda_\alpha=1$. In that case the
minimal scaled mse, attained when $J=\Lpot$, equals $N/S$, as it does
in the damped spherical harmonic case. We use this value as a lower
cutoff of the vertical axis on the left and the right.

\begin{figure*}\center
\rotatebox{0}{
\iftwocol
{\includegraphics[width=0.685\textwidth]{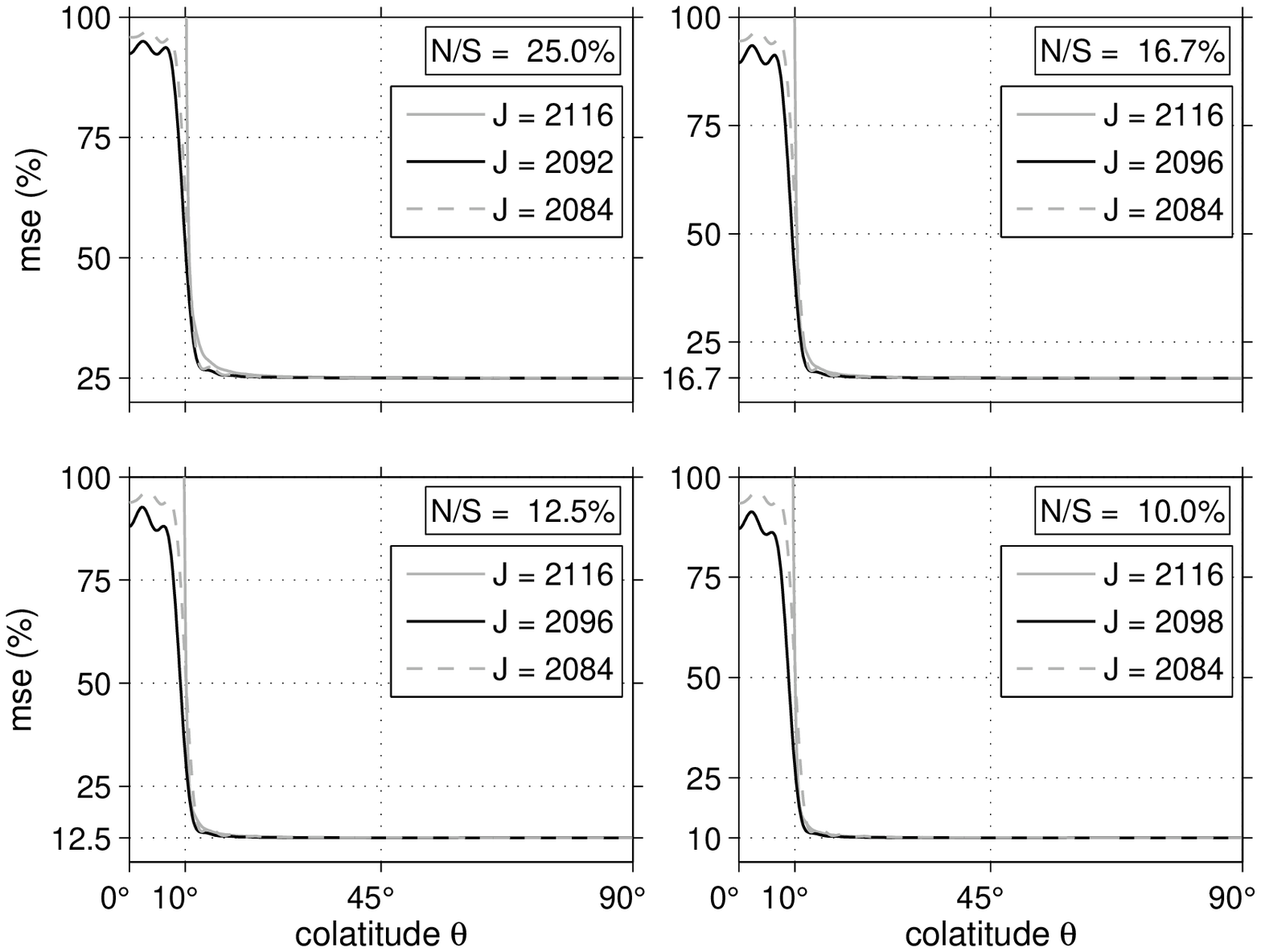}}
{\includegraphics[width=0.85\textwidth]{sderks.eps}}
} 
\caption{Mean-square error (as a percentage of the mean-square
signal strength) of the truncated Slepian function solution to the
geodetic estimation problem for a white signal and white noise. As in
Figure~\ref{sdkay}, the signal-to-noise levels shown are $S/N= 4, 6,
8$ and $10$, the measurement altitude $a=0$, the polar gap consists of
caps with radius $\Theta=10^{\circ}$, and the bandwidth of the signal
and its estimate is $L=45$. As a function of colatitude, for an
arbitrary longitude, we plot the mse of the untruncated solution, when
$J=\Lpot$, of truncation at the Shannon number, $J=K$, and of the
truncation levels which minimize the mse,
i.e. $J_\Omega=J_R=2092,2096,2096$, and $2098$, that were marked by
black vertical lines in Figure~\ref{sdkay}. The mse is symmetric about
the equator. The ordinate is truncated at 100\%; the mse of the
untruncated solution in the region of the polar gap exceeds this value
by several orders of magnitude.
}
\label{sderks}  
\end{figure*}
\begin{figure*}\center
\rotatebox{0}{
\iftwocol
{\includegraphics[width=0.685\textwidth]{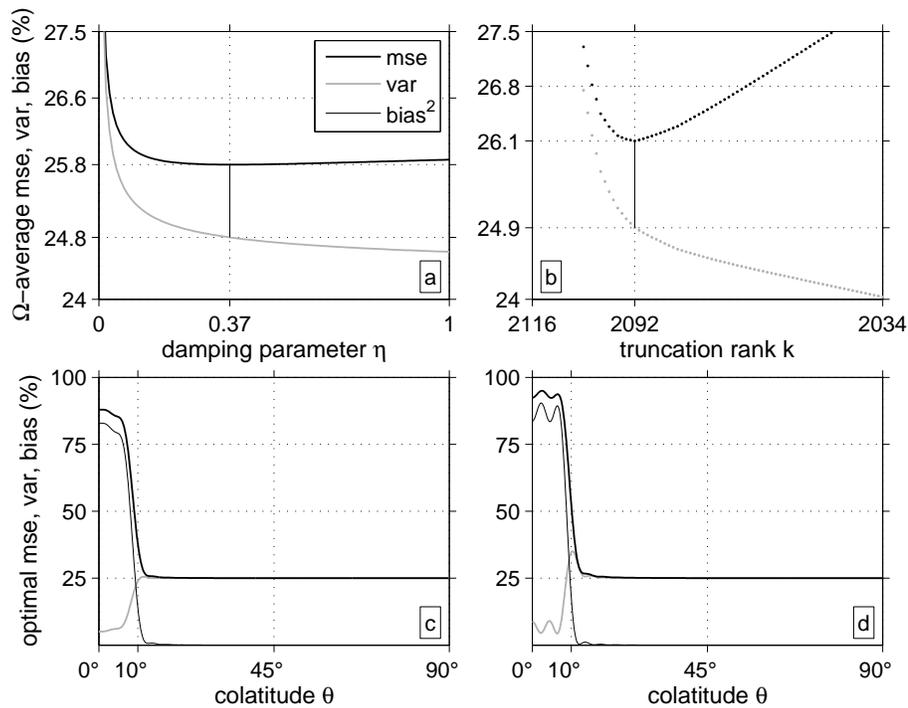}}
{\includegraphics[width=0.85\textwidth]{sdhoerl.eps}}
}
\caption{Mean square error (mse), variance and bias for the damped
least-squares solution and the truncated Slepian approach. The
antipodal polar gap has a radius $\Theta=10^{\circ}$; the bandwidth is
$L=45$; signal-to-noise ratio $S/N=4$. Panels a \& b show the values
averaged over the unit sphere, as in Figures~\ref{sdeta}
and~\ref{sdkay}; the squared bias is the difference between the mse
and variance curves, indicated by the thin black vertical line. Panels
c \& d show the values at the optimal damping and truncation
levels. The minimization of the mse reflects the trade-off between
variance, which dominates the uncovered regions of the polar gap, and
bias, dominant in the covered equatorial belt region.  }
\label{sdhoerl}  
\end{figure*} 

We may obtain the truncation level that minimizes the $\Omega$-average
mse by minimizing eq.~(\ref{SGmseO}) with respect to $J$. This will
yield an optimal truncation level $J_\Omega$. Likewise, minimizing
eq.~(\ref{SGmseR}) returns the truncation value $J_R$ at which the
$R$-average mse is minimal. Both minimization problems result in
identical constraints on the eigenvalue of the $J$-th eigenfunction
beyond which we truncate:
\be
\lambda_{J_\Omega}=\lambda_{J_R}\approx\frac{N}{S}
,
\ee
which is implicit but solvable. In Figure~\ref{sdkay}, the values
$J_\Omega=J_R$ identified by the top labels, are shown as a single
solid black vertical line; the Shannon number, $K=2048$, is shown by
the dotted black line and the bottom labels.

The mse~(\ref{SGmsefinala0}) for axisymmetric coverage regions is
independent of the longitude, as can be understood from
eq.~(\ref{polarg2}). Thus, in Figure~\ref{sderks} we plot
$\langle\epsilon^2(\theta)\rangle/\langle s^2(\theta)\rangle$, in
percent, as a function of colatitude for various truncation levels:
the untruncated, $J=\Lpot$, and optimally truncated cases,
$J=J_\Omega=J_R$, and the case truncated at the Shannon number
$J=K$. The vertical axes are truncated at 100\% since the untruncated
values exceed this value by many orders of magnitude. 

Finally, in Figure~\ref{sdhoerl}, we plot the relative contributions
of variance and bias for both the damped spherical harmonic and the
truncated Slepian case. The top panels show the mean square error, the
variance and the squared bias, according to the
relation~(\ref{nepsi}), as a function of the damping level (top left)
or the truncation rank (top right). The bottom panels show the
breakdown of mse, variance and bias in the spatial domain. For both
estimation methods, the bias is predominantly concentrated in the
areas over which no measurements are available, where it is generated
by the power of the ``missing signal''. The variance, on the other
hand, arises in the areas of coverage and is influenced by the power
of the noise. Both estimation methods show very similar results. Over
the covered area, the mse is nearly identical, and in the uncovered
region, the mse of the truncated Slepian case approaches that of the
damped spherical harmonic case, but it is slightly higher. It is
remarkable that eqs~(\ref{msefin2}) and~(\ref{SGmsefinala0}), despite
their different form, both give rise to a nearly complete spatial
separation of bias and variance. Only in eq.~(\ref{SGmsefinala0}) is
this separation immediately obvious by inspection: signal strength,
and thus bias, affect the low-ranking Slepian functions, whose power
is mostly concentrated inside of the polar gap, whereas noise, and
thus variance, affect the high-ranking Slepian functions whose power
is localized to the area of satellite coverage.

\section{C~o~n~c~l~u~s~i~o~n~s}
\label{c}

Spherical Slepian functions provide a natural solution to the geodetic
problem of having a polar gap in the satellite coverage of planetary
gravitational or magnetic fields. Indeed, the ill-posed geodetic
estimation problem of finding the source-level potential from noisy
observations taken at an altitude over an incomplete region of
coverage has natural connections to Slepian's spherical problem of
spatiospectral localization. We have proposed a new method that
expands the source field in terms of a truncated basis set of
spherical Slepian functions, and compared its statistical performance
with the traditional damped least-squares method in the spherical
harmonic basis. The optimally truncated Slepian method performs nearly
as well as the optimally damped spherical harmonic method, but it has
the significant advantage of an intuitive separation of the estimation
bias and variance over those Slepian functions sensitive to the
uncovered and covered regions, respectively. The construction of
Slepian functions over axisymmetric domains such as the latitudinal
belt or its complement, the polar gap, previously dismissed as
computationally unstable, has been shown to be eminently tractable. We
have shown that the operator that bandlimits a field on the unit
sphere and projects it onto the polar caps commutes with a
Sturm-Liouville operator. Its eigenfunctions, the Slepian functions,
can be computed extremely accurately and efficiently by diagonalizing
a tridiagonal matrix with analytically prescribed elements. The gains
in ease, speed, and accuracy thus achieved makes the use of spherical
Slepian functions in earth and planetary geodesy practical, as our
examples have shown.

\acknowledgments
We thank Mark Wieczorek for constructive coments on a preliminary
draft. This work was supported by a NERC Young Investigators' Award
and a Nuffield Foundation grant for Newly Appointed Lecturers to FJS,
and by Grant EAR-0105387 from the U.S.~National Science Foundation to
FAD.

\bibliography{/home/fjsimons/BIBLIO/bib.bib}
\bibliographystyle{gji}

\appendix
\section{Gr\"unbaum commutation}
\label{comm}

In this section we prove that the differential operator of
eq.~(\ref{grunbaumop2}), rewritten for $\mu=\cos\theta$ and
$b=\cos\Theta$, 
\ber
\ssT_p&=&(b^2-\mu^2)\nabla_m^2-2\mu(1-\mu^2)\frac{d}{d\mu}
-L_p(L_p+3)\mu^2\nnr\\ 
&=&\fracd{d}{d\mu}\bigg[(b^2-\mu^2)
(1-\mu^2)\fracd{d}{d\mu}\bigg]-L_p(L_p+3)\mu^2\nnr\\
&&{}-\fracd{m^2(b^2-\mu^2)}{1-\mu^2}
,
\label{Gproof}
\eer
commutes with the integral operator acting on $h_p(\mu')$
in eq.~(\ref{Fredholm2}),
\be
\int_{b}^1D_p(\mu,\mu')\,h_p(\mu')\,d\mu'=
\lambda\hspace{0.05em}h_p(\mu)
,
\label{Iproof}
\ee
whose symmetric kernel, $D_p(\mu,\mu')$, is given by
eq.~(\ref{Dproof1}), rewritten here for convenience in the form
\be
D_p(\mu,\mu')=\sumpp \Alma P_{lm}(\mu)P_{lm}(\mu')
.
\label{Dproof}
\ee
The domain of eq.~(\ref{Iproof}) is the interval~(\ref{thvalid}) of the double cap
\be
\{\mu: b \le\mu\le 1\}\cup\{\mu: -1 \le\mu\le -b\}
.
\ee
As in eqs~(\ref{plm})--(\ref{remind}), $P_{lm}$ is the associated
Legendre polynomial of degree $l$ and order $m$, and $\Alma$ is a
normalization constant. We remind the reader of our notation:
$h_p(\mu)$ is a colatitudinally dependent function that is limited in
space to the antipodal polar caps of radius $\Theta=\cos^{-1} b$. It
is either odd or even about the equator, as indicated by the subscript
$p$:
\begin{subequations}
\ber
h_p(\mu)&=&0, \quad-b\ge\mu\ge b
,\\
h_e(\mu)&=&h_e(-\mu)
,\\
h_o(\mu)&=&-h_o(-\mu)
.
\eer
\end{subequations}
The solutions to eq.~(\ref{Iproof}) are functions $h_p(\mu)$ that are
spectrally concentrated in a spherical harmonic degree interval $0\le
l\le L$; the eigenvalue $\lambda$ is the quadratic measure of this
concentration~(\ref{normratio2}). The primed summation symbol skips
every other term in the interval from $m_p$ to $L_p$, which are both
of the same parity, either even or odd. Depending on the requested
order $m$ and concentration bandwidth $L$ of the solutions, $m_p$ is
either $m$ or $m+1$, and $L_p$ is either $L$ or $L-1$, following
eq.~(\ref{Lpdef}). We further distinguish $\ssT$ acting on $\mu$ from
$\ssT'$ that acts on $\mu'$.

To confirm commutativity we are required to show that
\ber
\lefteqn{\int_{b}^1 D_p(\mu,\mu') \,\ssT'_p h_p(\mu')\,d\mu'=}\nnr\\
&&{}\hspace{6em}
\int_{b}^1 \ssT_p D_p(\mu,\mu')\,h_p(\mu')\,d\mu' 
.
\label{grunshow0}
\eer
We first show that the left side of eq.~(\ref{grunshow0}) can be
rewritten as 
\ber
\lefteqn{\int_{b}^1 D_p(\mu,\mu') \,\ssT'_p h_p(\mu')\,d\mu'=}\nnr\\
&&{}\hspace{6em}
\int_{b}^1 \ssT'_p D_p(\mu,\mu')\,h_p(\mu')\,d\mu',
\label{grunshow1}
\eer
and then we verify that
\be
\ssT_pD_p(\mu,\mu')=\ssT'_pD_p(\mu,\mu').
\label{grunshow2}
\ee
The first result~(\ref{grunshow1}) is easily verified by integration
by parts: for any two functions $\zeta(\mu)$ and $\eta(\mu)$, it may
be shown that, whether $b=\cos\Theta$ or $b=-1$,
\ber
\int_{b}^1\zeta\,(\ssT_p\eta)\,d\mu&=&
-\int_{b}^1\Big[(b^2-\mu^2)(1-\mu^2)
\frac{d\zeta}{d\mu}\frac{d\eta}{d\mu}\nnr\\
&&\quad{}+L_p(L_p+3)\mu^2\,\zeta\eta \nnr\\
&&\quad{}+\,m^2(b^2-\mu^2)(1-\mu^2)^{-1}
\zeta\eta\Big]\,d\mu\nnr\\ 
&=&\int_{b}^1(\ssT_p\hsp\zeta)\,\eta\,d\mu.
\label{grunshow3}
\eer
To verify the second result~(\ref{grunshow2}) we use the
Laplace-Beltrami identity $\nabla_{\!m}^2P_{lm}=-l(l+1)P_{lm}$
\cite[]{Dahlen+98} to write
\iftwocol
{
\ber
\label{grunshow4}
\lefteqn{(\ssT_p-\ssT'_p)D_p(\mu,\mu')=}\hspace{1em}\\
&&{}(\mu^2-\mu'^2)\sumpp
\Alma P_{lm}(\mu)P_{lm}(\mu')\nnr\\
&&{}\hspace{6.5em}
\times\big[l(l+1)-L_p(L_p+3)\big]\nnr \\
&&{}-2\mu(1-\mu^2)\sumpp\Alma
\frac{d}{d\mu}P_{lm}(\mu)P_{lm}(\mu')\nnr \\
&&{}+2\mu'(1-\mu'^2)\sumpp\Alma 
P_{lm}(\mu)\frac{d}{d\mu'}P_{lm}(\mu').\nnr
\eer
}
{
\ber
\label{grunshow4}
\lefteqn{(\ssT_p-\ssT'_p)D_p(\mu,\mu')=}\hspace{1em}\\
&&{}(\mu^2-\mu'^2)\sumpp
\Alma P_{lm}(\mu)P_{lm}(\mu')
\big[l(l+1)-L_p(L_p+3)\big]\nnr \\
&&{}-2\mu(1-\mu^2)\sumpp\Alma
\frac{d}{d\mu}P_{lm}(\mu)P_{lm}(\mu')
+2\mu'(1-\mu'^2)\sumpp\Alma 
P_{lm}(\mu)\frac{d}{d\mu'}P_{lm}(\mu').\nnr
\eer
}
Two well-known Legendre identities help us simplify the above; a
derivative identity and a recursion relation \cite[]{Dahlen+98}: 
\begin{subequations}
\ber
\label{grunshow5}
\frac{dP_{lm}}{d\mu}&=&\frac{(l+1)\mu P_{lm}-(l-m+1)P_{l+1\,m}}{1-\mu^2}\\
\mu P_{lm}&=&\frac{(l-m+1)P_{l+1\,m}+(l+m)P_{l-1\,m}}{2l+1}
\label{latest}.
\eer
\end{subequations}
Applying eq.~(\ref{grunshow5}), then eq.~(\ref{latest}),
transforms eq.~(\ref{grunshow4}) into
\ber
\label{grunshow6}
\lefteqn{(\ssT_p-\ssT'_p)D_p(\mu,\mu')=}\\
&&{}(\mu^2-\mu'^2)\sumpp
\Alma P_{lm}(\mu)P_{lm}(\mu')\nnr\\
&&{}\hspace{3em}
\times\big[(l-2)(l+1)-L_p(L_p+3)\big]\nnr \\
&&{}+\sumpp\Alma\big[P_{l+2\,m}(\mu)P_{lm}(\mu')-P_{lm}(\mu)
P_{l+2\,m}(\mu')\big]\nnr\\
&&{}\hspace{3em}
\times2\,\frac{(l-m+1)(l-m+2)}{2l+3}\nnr
.
\eer
The Legendre derivative identity of eq.~(\ref{grunshow5}) can be
manipulated to yield a formula of the Christoffel-Darboux type
\cite[]{Simons+2006a}:
\ber 
\lefteqn{
(\mu^2-\mu'^2)\sumpp\Alma P_{lm}(\mu)P_{lm}(\mu')=}\label{christdarboux}\\
&&{}\big[P_{L_p+2\,m}(\mu)P_{L_pm}(\mu')-P_{L_pm}(\mu)
P_{L_p+2\,m}(\mu')\big]\nnr\\
&&{}\hspace{3em}
\times\fracd{(L_p-m+2)!}{(2L_p+3)(L_p+m)!}\nnr
.
\eer
Inserting eq.~(\ref{christdarboux}) into eq.~(\ref{grunshow6}) yields:
\ber
\lefteqn{(\ssT_p-\ssT'_p)D_p(\mu,\mu')=}\label{grunshow7}
\\
&&{}\hspace{-1em}(\mu^2-\mu'^2)\sumpp\Alma P_{lm}(\mu)P_{lm}(\mu')\nnr\\
&&{}\hspace{-1em}\times\big[(l-2)(l+1)-L_p(L_p+3)\big]\nnr\\
&&{}\hspace{-1em}+(\mu^2-\mu'^2)\sumpp2(2l+1)\sumpn\Almb
P_{nm}(\mu)P_{nm}(\mu')\nnr
.
\eer
Interchanging the order of summation and relabeling the sums,
\ber
\lefteqn{(\ssT_p-\ssT'_p)D_p(\mu,\mu')=}\label{grunshow8}\\
&&{}(\mu^2-\mu'^2)\sumpp\Alma P_{lm}(\mu)P_{lm}(\mu')\nnr\\
&&{}\times\left[(l-2)(l+1)-L_p(L_p+3)+\sump_{n=l}^{L_p}2(2n+1)\right]\nnr
.
\eer
The term in square brackets always vanishes; in other words,
$(\ssT_p-\ssT'_p)D_p(\mu,\mu')=0$ and the commutation relation of
eq.~(\ref{grunshow0}) is confirmed. Since commuting operators have the
same eigenfunctions, we can find the spacelimited, fixed-order
eigenfunctions $h(\theta)$ or $h(\mu)$ by solving the integral
equation eq.~(\ref{Iproof}), or, equivalently, by solving the
differential equation
\be
\ssT_ph_p(\mu)=\chi\hspace{0.05em}h_p(\mu)
\label{Gchih}
,
\ee
on the domain of the double polar cap, where $\chi\not= \lambda$ is
the associated eigenvalue. The operator $\ssT_p$ is Sturm-Liouville,
i.e., it is of the form
\be
\label{sturm}
(ph')'-qh+\chi\rho h =0,
\quad \cos\Theta\le\mu\le 1,
\ee
where $p(\mu)=(\mu^2-b^2)(1-\mu^2)$,
$q(\mu)=m^2(1-\mu^2)^{-1}(\mu^2-b^2)-L_p(L_p+3)\mu$, $\rho(\mu)=1$, 
and the prime denotes differentiation with respect to $\mu$. Since
$\ssT_p$ is a Sturm-Liouville operator, the eigenvalue spectrum of
eq.~(\ref{Gchih}) is simple: the spacing between adjacent fixed-order
eigenvalues is roughly equant, as we illustrated in
Figure~\ref{sdgrunval}. 

\section{The Gr\"unbaum matrix}

The domains of eqs~(\ref{Iproof}) or~(\ref{Gchih}), originally
restricted to the area contained within the double polar caps, may be
extended to the entire sphere, $0\le \cos\theta\le\pi$, by writing
\be
\ssT_pg_p(\mu)=\chi\hspace{0.05em}g_p(\mu)\quad,\quad -1\le\mu \le1
\label{Gchig}
.
\ee
The unknown functions, see eq.~(\ref{grep}), must now be bandlimited rather than
spacelimited:
\be
g_p=\sumpp g_{lm}\Xlmth
\label{grep2}
,
\ee
where, as in eq.~(\ref{Lpdef}), depending on the requested order $m$
and concentration bandwidth $L$ of the solutions, $m_p$ is either $m$
or $m+1$, $L_p$ is either $L$ or $L-1$, and the primed summation skips
every second integer. As a result, $g_{lm}$ requires no further
identification. Inserting the representation of eq.~(\ref{grep2}) into
eq.~(\ref{Gchig}), multiplying both sides by $2\pi\sin\theta\,\Xlmth$,
integrating over all colatitudes $0\le\theta\le\pi$, and invoking the
orthogonality eq.~(\ref{Xlmortho}) easily transforms eq.~(\ref{Gchig})
into an algebraic eigenvalue equation for the coefficients of the
unknown functions $g$:
\be
\sfT_p\sfg_p=\chi\sfg_p
,
\label{Gmatrix}
\ee
where we define
\be
T^p_{ll'}=2\pi\int_{0}^{\pi}\Xlm(\ssT_p\hspace{-0.05em}
\Xlpm)\sin\theta\,d\theta
.
\label{Gmatrix2}
\ee
To avoid clutter, we have changed the subscript $p$ indicating the
parity of the solutions into a superscripted $p$ on the matrix
elements $T^p_{ll'}$. These are indexed by the integer degrees $l$ and
$l'$. Since in $\sfT_p$ the only degrees involved range from $m_p$ to
$L_p$, with every second degree skipped, when $m=0$, the first element
of the matrix $\sfT_e$ is thus $T^e_{00}$, the second $T^e_{20}$, and
so on, whereas the first and second elements of $\sfT_o$ are
$T^o_{11}$ and $T^o_{30}$, and so on, respectively. The matrix whose
diagonalization leads to the even functions $g_e$ is thus given by
\be
\sfT_e =\left(\begin{array}{cccc}
T^e_{mm} & T^e_{mm+2} & \cdots & T^e_{mL\!_e}\\
T^e_{m\!+2\,m} & T^e_{m\!+2\,m\!+2} & \cdots & T^e_{m\!+2\,L\!_e}\\
\vdots & \vdots & \vdots\\
T^e_{L\!_em} & T^e_{L\!_e\,m\!+2} & \cdots & T^e_{L\!_eL\!_e}\end{array}\right)
.
\label{Ge}
\ee
Its dimensions are $[(L_e-m)/2+1]\times [(L_e-m)/2+1]$.  The
 $[(L_o-m-1)/2+1]\times [(L_o-m-1)/2+1]$--dimensional matrix yielding
 the odd functions $g_o$ is
\be
\sfT_o =\left(\begin{array}{cccc}
T^o_{m\!+1m\!+1} & T^o_{m\!+1m+3} & \cdots & T^o_{m\!+1L\!_o}\\
T^o_{m\!+3\,m\!+1} & T^o_{m\!+3\,m\!+3} & \cdots & T^o_{m\!+3\,L\!_o}\\
\vdots & \vdots & \vdots\\
T^o_{L\!_om\!+1} & T^o_{L\!_o\,m\!+3} & \cdots &
T^o_{L\!_oL\!_o}\end{array}\right).
\label{Go}
\ee
The dimensions of the combined matrix 
\be
\sfT'={\rm diag}\,\left(\sfT_e,\sfT_o\right)
\ee
are thus $(L-m+1)\times (L-m+1)$, as expected. Such a matrix can be
thought of as a permutation of a ``full form'' 
\be
\sfT= \left(\begin{array}{ccc}
T_{mm} & \cdots & T_{mL}\\
\vdots && \vdots\\
T_{Lm} & \cdots & T_{LL}\end{array}\right)
,
\ee
where the degrees are arranged in the correct order, without skipping
entries. It is this matrix $\sfT$ that commutes with the matrix
$\sfD$ of eq.~(\ref{kernel4}), which is also in the form
\be
\sfD = \left(\begin{array}{ccc}
D_{mm} & \cdots & D_{mL}\\
\vdots && \vdots\\
D_{Lm} & \cdots & D_{LL}\end{array}\right)
\ee
and whose elements are given by eq.~(\ref{kernel4}),
\be
D_{ll'}=2\pi\int_{0}^{\Theta}\Xlm\Xlpm\sin\theta\,d\theta.
\label{kernelv}
\ee
Deriving the form of the matrix entries $T^p_{ll'}$ of
eq.~(\ref{Gmatrix2}) requires the  Gr\"unbaum operator $\ssT_p$  of
eq.~(\ref{Gproof}) as a function of colatitude, i.e.
\ber
\ssT_p&=&(b^2-\cos^2\theta)\nabla_m^2+
2\cos\theta\sin\theta\frac{d}{d\theta}\nnr\\ 
&&{}-L_p(L_p+3)\cos^2\theta
.
\eer
Evaluating eq.~(\ref{Gmatrix2}) is perhaps not as pedestrian as by
``simply reading off directly the inner products''  as claimed by
\cite{Grunbaum+82}, but since only the result is of any practical
consequence here, we will simply state it:
\begin{subequations}\ber
T^p_{ll}&=&-l(l+1)\cos^2{\Theta}
+\frac{2}{2l+3}\left[(l+1)^2-m^2\right]\nnr\\
&&{}+[(l-2)(l+1)-L_p(L_p+3)]\nnr\\
&&{}\times\left[\frac{1}{3}-\frac{2}{3}\,
\fracd{3m^2-l(l+1)}{(2l+3)(2l-1)}\right],\\ 
T^p_{l\,l+2}&=&\fracd{\big[l(l+3)-L_p(L_p+3)\big]}{2l+3}\\
&&{}\times\sqrt{\fracd{\left[(l+2)^2-m^2\right]
\left[(l+1)^2-m^2\right]}{(2l+5)(2l+1)}},\nnr\\
T^p_{ll'}&=&0\quad\mbox{otherwise}
.
\eer\end{subequations}
This is the result given as eq.~(\ref{gdefi}) in the main text. The
elements again are specified in terms of the degree and not by the
matrix index. Both $\sfT_e$ and $\sfT_o$ are thus not only real and
symmetric, but also tridiagonal. The important result is that the
coefficients of the even or odd optimally concentrated antipodal polar
cap eigenfunctions $g_e(\theta)$ or $g_o(\theta)$ both only require
the numerical diagonalization of a tridiagonal matrix. Both matrices
have analytically prescribed elements, and a spectrum of eigenvalues
that is guaranteed to be simple.

\label{lastpage}
\end{document}

\end{document}